\documentclass{article}

\usepackage{amsfonts}
\usepackage{amssymb}
\usepackage{textcomp}
\usepackage{latexsym,amsmath}
\usepackage{graphics}

\usepackage{amsfonts}
\usepackage{amssymb}
\usepackage{textcomp}
\usepackage{latexsym,amsmath}

\renewcommand{\phi}{\varphi}
\newcommand{\be}{\begin{equation}}
\newcommand{\ee}{\end{equation}}
\newcommand{\ba}{\begin{eqnarray}}
\newcommand{\ea}{\end{eqnarray}}
\newcommand{\ban}{\begin{eqnarray*}}
\newcommand{\ean}{\end{eqnarray*}}

\newcommand{\nul}{{\bf0}}
\newcommand{\rd}{{\mathbb R}^d}
\newcommand{\zd}{{\mathbb Z}^d}
\newcommand{\td}{{\mathbb T}^d}

\newcommand{\z} {{\mathbb Z}}

\newcommand{\n} {{\mathbb N}}

\newcommand{\h}{\widehat}
\newcommand{\w}{\widetilde}
\newcommand{\too}{\mathop{\longrightarrow}}

\def\supp{\operatorname{supp}}
\def\sinc{\operatorname{sinc}}
\def\const{\operatorname{const}}

\def\N{{{\Bbb N}}}
\def\Z{{{\Bbb Z}}}
\def\T{{{\Bbb T}}}
\def\R{{\Bbb R}}

\def\vp{{\varphi}}

\def\({\left(}
\def\){\right)}

\def\C{{\Bbb C}}
\def\l{{\lambda }}
\def\a{{\alpha }}
\def\D{{\Delta }}

\def\a{{\alpha}}
\def\b{{\beta}}
\def\d{{\delta}}

\def\vp{{\varphi}}

\def\g{{\gamma }}

\def\L{{\Lambda}}
\def\F{{\mathcal{F}}}
\def\A{{\mathcal{A}}}
\def\S{{\mathcal{S}}}

\newtheorem{theo}{Theorem}
\newtheorem{lem}[theo]{Lemma}

\newtheorem {coro} [theo] {Corollary}
\newtheorem {defi} [theo] {Definition}
\newtheorem {rem} [theo] {Remark}

\title{Approximation by  multivariate quasi-projection operators \\  and Fourier multipliers
}
\author{
Yurii Kolomoitsev$^{1, 2}$ and Maria Skopina$^{3, 4}$
}
\date{\small $^{1}$Universit\"at zu L\"ubeck,
L\"ubeck, Germany; kolomoitsev@math.uni-luebeck.de  \\
\small $^{2}$Institute of Applied Mathematics and Mechanics of NAS of Ukraine, Slov'yans'k, Ukraine
\\
$^{3}$Saint Petersburg State University,  St. Petersburg,  Russia;
skopina@ms1167.spb.edu \\
$^{4}$Regional Mathematical Center of Southern Federal University
}

\begin{document}

\maketitle

\begin{abstract}
Multivariate quasi-projection operators  $Q_j(f,\varphi, \widetilde{\varphi})$, associated with a function $\phi$ and a distribution/function $\widetilde{\varphi}$, are considered. The function $\varphi$ is supposed to satisfy  the Strang-Fix conditions  and a  compatibility condition with $\widetilde{\varphi}$. Using technique based on the Fourier multipliers, we stu\-died approximation properties of such operators  for functions  $f$ from anisotropic Besov spaces and $L_p$ spaces with $1\le p\le \infty$. In particular, upper and lower estimates of the $L_p$-error of approximation  in terms of anisotropic moduli of smoothness and  anisotropic best approximations are obtained.
\end{abstract}

\bigskip


\medskip

	
\section{Introduction}

The multivariate quasi-projection operator with a matrix dilation  $M$  is defined as:
$$%\be
 Q_j(f, \phi,\w\phi)=|\det M|^{j}  \sum_{n\in\zd} \langle f,\w\vp(M^j\cdot+n)\rangle \vp(M^j\cdot+n),
%\label{my00}
$$%\ee
%???where $\phi$ is a function and $\w\phi$ is a tempered distribution/function,  and the functional
%$\langle f,\w\vp(M^j\cdot+n) \rangle$
%has meaning in some sense. ???
where $\phi$ is a function, $\w\phi$ is a tempered distribution, and 
$\langle f,\w\vp(M^j\cdot+n) \rangle$
is an appropriate functional.

The class of  operators $Q_j(f, \phi,\w\phi)$ is quite large. It includes the operators associated with a regular function $\w\phi$,
in particular, the so-called scaling expansions appearing in wavelet constructions
(see, e.g.,~\cite{BDR, v58, Jia2, KPS, KS, Sk1}) as well as the
Kantorovich-Kotelnikov operators  and their generalizations (see, e.g.,~\cite{CV2, KS3, KS20, OT15, VZ2}).
An essentially different class consists of the operators $Q_j(f,\phi,\w\phi)$ associated with a tempered distribution $\w\phi$  related to the Dirac delta-function (the so-called sampling-type operators). The model example of such operators is the following classical sampling expansion, appeared originally in the Kotelnikov formula,
$$
\sum_{n\in\z}  f(-2^{-j}n)\,\frac{\sin\pi(2^jx+n)}{\pi (2^jx+k)}=
2^{j}  \sum_{n\in\z} \langle f,\delta(2^j\cdot+n) \rangle\, {\rm sinc}(2^jx+n),
$$
where  $\delta$ is the Dirac delta-function and ${\rm sinc}\,x:=\frac{\sin\pi x}{\pi x}$.
In recent years, many authors have studied  %*** A lot of authors studied  
approximation properties of  the sampling-type operators for various functions $\phi$
(see, e.g.,~\cite{BD, Butz5, KM, KKS, KS20, KS, SS, Unser}).
Consideration of  functions $\phi$ with a good decay is very useful for different engineering applications. In particular, the operators associated with  a  linear combination of $B$-splines as  $\phi$, and  the Dirac delta-function as $\w\phi$, was studied, e.g.,
in~\cite{Butz6, Butz7, SS}. For a class of fast decaying functions $\phi$, the sampling-type
quasi-projection operators were considered in~\cite{KS}, where the error estimates in the $L_p$-norm,
$p\ge2$,  were obtained in terms of the Fourier transform of $f$, and the approximation
order of the operators was found in the case of an isotropic matrix $M$. These results were extended
to an essentially wider class of functions $\phi$ in~\cite{cksv} (see Theorem~A below).  Next,  in the paper~\cite{KS19}, the results of~\cite{KS} were improved  in several directions. Namely,  the  error estimates were obtained also for the case $1\le p<2$,
the requirement on the approximated function $f$ were weakened, and the estimates were
given in terms  of anisotropic moduli of smoothness and best approximations.

The main goal of the present paper is to extend the results of~\cite{KS19} to band-limited functions $\phi$ and to the case $p=\infty$. The scheme of the proofs of our results is similar to the one given in~\cite{KS19}, but the technic is essentially refined by means of using Fourier multipliers. This development allows also to improve the results for the class of fast decaying functions $\phi$ and to obtain lower estimates for the $L_p$-error of approximation by quasi-projection operators in some special cases. Similarly, the main result of~\cite{KS3} (see Theorem~B below) is essentially extended in several directions (lower estimates, fractional smoothness, approximation in the uniform metric).

The paper is organized as follows. Notation  and preliminary information are given in Sections~2 and 3, respectively.
Section~4 contains auxiliary results.  The main results are presented in Section~5. In particular, the $L_p$-error estimates for quasi-projection operators $Q_j(f,\phi,\w\phi)$ in the case of weak compatibility of $\phi$ and $\w\phi$ are obtained in Subsection~5.2. In this subsection, we also consider lower estimates for the $L_p$-error and a generalization of compatibility conditions to the case of fractional smoothness. Subsection~5.3 is devoted to approximation by operators $Q_j(f,\phi,\w\phi)$ in the case of strict compatibility $\phi$ and $\w\phi$. Two generalizations of the Whittaker--Nyquist--Kotelnikov--Shannon-type theorem are also proved in this subsection.

\section{Notation}

%We use the standard multi-index notations.
    As usual, we denote by $\n$  the set of positive integers, $\rd$ is the $d$-dimensional Euclidean space,
    $\zd$ is the integer lattice  in $\rd$, $\z_+^d:=\{x\in\zd:~x_k\geq~{0}, k=1,\dots, d\}$, and
    $\td=\rd\slash\zd$ is the $d$-dimensional torus.
    Let  $x = (x_1,\dots, x_d)^{T}$ and
    $\xi =(\xi_1,\dots, \xi_d)^{T}$ be column vectors in $\rd$,
    then $(x, \xi):=x_1\xi_1+\dots+x_d\xi_d$,
    $|x| := \sqrt {(x, x)}$,  $\nul=(0,\dots, 0)^T\in \rd$,
		and $B_\d=\{x\in \R^d\,:\,|x|< \d\}$.
		%denotes the ball of radius $\d>0$ with the center in $\nul$.

Given $a,b\in\rd$ and $\alpha\in\zd_+$, we set
    $$
    [\alpha]=\sum\limits_{k=1}^d \alpha_k, \quad
    D^{\alpha}f=\frac{\partial^{[\alpha]} f}{\partial x^{\alpha}}=\frac{\partial^{[\alpha]} f}{\partial^{\alpha_1}x_1\dots
    \partial^{\alpha_d}x_d}, \,\,\quad
 a^b=\prod\limits_{j=1}^d a_j^{b_j},\,\,\quad
  \alpha!=\prod\limits_{j=1}^d\alpha_j!\,.
    $$

If $M$ is a $d\times d$ matrix,
then $\|M\|$ denotes its operator norm in $\rd$; $M^*$ denotes the conjugate matrix to $M$, $m=|\det M|$.
By $I$ we denote the identity matrix, i.e., $I=M^0$.

A $d\times d$ matrix $M$ is called a  dilation matrix if its eigenvalues are bigger than one in modulus.
We denote the set of all dilation matrices by $\mathfrak{M}$.
It is well known that $\lim_{j\to\infty}\|M^{-j}\|=0$ for dilation matrices.
For any $M\in \mathfrak{M}$, we set $\mu_0:=\min\{\mu\in \N\,:\,\T^d\subset \frac12 M^{*\nu}\T^d\quad\text{for all}\quad\nu\ge \mu-1\}$.

Recall that a  matrix $M$ is  isotropic if it is similar to a multiple of an orthogonal matrix,  its  eigenvalues  $\lambda_1,\dots,\lambda_d$ are such that $|\lambda_1|=\cdots=|\lambda_d|$.

As usual, $L_p$ denotes the space $L_p(\rd)$, $1\le p\le \infty$, with the norm $\Vert \cdot\Vert_p=\Vert \cdot\Vert_{L_p(\R^d)}$,
${C}$ denotes the space of all uniformly continuous bounded functions on $\R^d$, and
$$
{C}_0:=\{f\in {C}\,:\, \lim_{|x|\to\infty}f(x)=0\}.
$$
We use $W_p^n$, $1\le p\le\infty$, $n\in\n$, to denote the Sobolev space on~$\rd$, i.e. the set of
functions whose derivatives up to order $n$ are in $L_p$, with usual Sobolev norm.

If $f$ and $g$ are functions defined on $\rd$ and $f\overline g\in L_1$,
then
$$
\langle  f, g\rangle:=\int_{\rd}f(x)\overline{g(x)}dx.
$$
The convolution of functions $f$ and $g$ is defined by
$$
f*g(x)=\int_{\R^d} f(t)g(x-t)dt.
$$
The Fourier transform of $f\in L_1$ is given by
$$
\mathcal{F}f(\xi)=\widehat
f(\xi)=\int_{\rd} f(x)e^{-2\pi i
(x,\xi)}\,dx.
$$
For any function $f$, we denote $f^-(x)=\overline{f(-x)}$ .

The Schwartz class of functions defined on $\rd$ is denoted by $\mathcal{S}$.
The dual space of $\mathcal{S}$ is $\mathcal{S}'$, i.e. $\mathcal{S}'$ is the space of tempered distributions.
Suppose $f\in \mathcal{S}$ and $\phi \in \mathcal{S}'$, then
$\langle f, \phi \rangle:=\phi(f)$.
For any  $\phi\in \mathcal{S}'$, we define $\overline{\vp}$ and $\vp^-$  by $\langle f, \overline{\phi}\rangle:=\overline{\langle \overline{f},  \phi\rangle}$, $f\in \mathcal{S}$, and $\langle f, {\phi}^-\rangle:=\overline{\langle {f}^-,  \phi\rangle}$, $f\in \mathcal{S}$, respectively. The Fourier transform of $\phi$ is
defined by $\langle \h f, \h \phi\rangle=\langle f, \phi\rangle$,
$f\in \mathcal{S}$. The convolution of $\phi \in \mathcal{S}'$ and $f\in \mathcal{S}$ is given by $f*\vp(x)=\langle f, \overline{\vp(x-\cdot)}\rangle=\langle f, \vp(\cdot-x)^-\rangle$.
For suitable functions/distributions $f$ and $h$,  we denote by $\Lambda_h (f)$ the following multiplier operator:
$$
\Lambda_h (f) := \mathcal{F}^{-1}(h \h f ).
$$
Next, for a fixed matrix $M\in \mathfrak{M}$ and  a function $\phi$, we define $\vp_{jk}$ by
$$
\phi_{jk}(x):=m^{j/2}\phi(M^jx+k),\quad j\in\z,\quad k\in\rd.
$$
For $\w\phi\in \mathcal{S}'$, $j\in\z$, and $k\in\zd$, we define $\w\phi_{jk}$ by
        $$
        \langle f, \w\phi_{jk}\rangle:=
        \langle f_{-j,-M^{-j}k},\w\phi\rangle,\quad  f\in \mathcal{S}.
        $$

By $\mathcal{S}_N'$, $N\ge 0$, we denote the set of all tempered distribution $\w\phi$
	whose Fourier transform $\h{\w\phi}$ is a measurable function on $\rd$
	such that $|\h{\w\phi}(\xi)|\le c({\w\phi})(1+ |\xi|)^{N}$
	 for almost all $\xi\in\rd$.

%By ${\cal L}_p$, $1\le p \le \infty$, we denote the set
Let $1\le p \le \infty$. We set
	$$
	{\cal L}_p:=
	\left\{
	\phi\in L_p\,:\, \|\phi\|_{{\cal L}_p}:=
	\bigg\|\sum_{l\in\zd} \left|\phi(\cdot+l)\right|\bigg\|_{L_p(\td)}<\infty
	\right\}.
	$$
It is not difficult to see that ${\cal L}_1=L_1,$ $\|\phi\|_p\le \|\phi\|_{{\cal L}_p}$, and
$\|\phi\|_{{\cal L}_q}\le \|\phi\|_{{\cal L}_p}$ for $1\le q \le p \le\infty$.

For any $d\times d$ matrix $A$, we introduce the space
$$
\mathcal{B}_{A,p}:=\{g\in L_p\,:\, \supp \h g\subset A^*\T^d \}
$$
and the corresponding anisotropic best approximations
$$
E_{A} (f)_p:=\inf\{\Vert f-P\Vert_p\,:\, P\in \mathcal{B}_{A,p}\}.
$$

Let $\alpha$ be a positive function defined on the set of all $d\times d$ matrices $A$.
We consider the following anisotropic Besov-type space associated with a matrix~$A$.
  We  say that
$f\in \mathbb{B}_{p;A}^{\a(\cdot)}$,  $1\le p\le\infty$, if $f\in L_p$ for $p<\infty$,  $f\in {C}_0$ for $p=\infty$, and
$$
\Vert f\Vert_{\mathbb{B}_{p;A}^{\a(\cdot)}}:=\Vert f\Vert_p+\sum_{\nu=1}^\infty |\det A|^\frac\nu p\a(A^\nu) E_{A^\nu} (f)_p<\infty.
$$
Note that in the case $A=2 I$ and $\a(\cdot)\equiv \a_0\in \R$, the space $\mathbb{B}_{p;A}^{\a(\cdot)}$ coincides with the classical Besov space $B_{p,1}^{d/p+\a_0}(\R^d)$.

For any matrix $M\in \mathfrak{M}$, we denote by  $\mathcal{A}_M$ the set of all positive functions $\a\,:\,\R^{d\times d}\to \R_+$ that satisfy the condition $\a(M^{\mu+1})\le c(M)\a(M^\mu)$ for all $\mu\in \Z_+$.

For any $d\times d$ matrix $A$, we introduce the anisotropic fractional modulus of smoothness of order $s$, $s>0$,
$$
\Omega_s(f,A)_p:=\sup_{|A^{-1}h|<1} \Vert \Delta_h^s f\Vert_p,
$$
where
$$
\Delta_h^s f(x):=\sum_{\nu=0}^\infty (-1)^\nu \binom{s}{\nu} f(x+h\nu).
$$
Recall that the standard fractional modulus of smoothness of order $s$, $s>0$, is defined
by
$$
\omega_s(f,t)_p:=\sup_{|h|<t} \Vert \Delta_t^s f\Vert_p,\quad t>0.
$$
We refer to~\cite{KT20} for the collection of basic properties of moduli of smoothness in $L_p(\R^d)$.

For an appropriate function $f$ and  $s>0$, the fractional power of Laplacian is given by
$$
(-\Delta)^{s/2} f(x):= \mathcal{F}^{-1}\(|\xi|^s\h f(\xi)\)(x).
$$

As usual, $\ell_p$, $1\le p\le\infty$, denotes the space of all sequences $a=\{a_n\}_{n\in \Z^d}\subset\C$ equipped with the norm
$$
\Vert a\Vert_{\ell_{p}}:=\left\{
                                                 \begin{array}{ll}
                                                  \displaystyle\bigg(\sum\limits_{n\in \Z^d}|a_n|^p\bigg)^\frac1p, & \hbox{if $p<\infty$,} \\
                                                   \displaystyle\sup_{n\in \Z^d}|a_n|, & \hbox{if $p=\infty$,}
                                                 \end{array}
                                               \right.
$$
and ${\rm c}_0$ denotes the subspace of  $\ell_\infty$ consisting of the sequnces converging to zero.

By $\eta$ we denote a real-valued function in $C^\infty(\R^d)$ such that $\eta(\xi)=1$ for $\xi \in \T^d$ and $\eta(\xi)=0$ for $\xi\not\in 2\T^d$. For any $\d>0$, we denote $\eta_\d=\eta(\d^{-1}\cdot)$.

Finally, for any  $p\in [1,\infty]$, we define $p'$ by $1/p'+1/p=1$ and write $c, c_1, c_2,\dots$ to denote positive constants that depend on indicated parameters.

\section{ Preliminary information and main definitions.}
\label{sa}

In what follows, we discuss the operator
$$
Q_j(f,\phi,\w\phi):=\sum_{k\in\zd} \langle f, {\w\phi}_{jk}\rangle \phi_{jk},
$$
where the "inner product"
$\langle f, {\w\phi}_{jk}\rangle$
has meaning in some sense. This operator is associated with a  matrix $M$, which is a  matrix dilation by default.

The operator $Q_j(f,\phi,\w\phi)$ is an element of the shift-invariant space generated by the function~$\phi$.
It is  known that a function $f$ may be approximated by the  elements of such shift-invariant space  only if $\phi$ satisfies the so-called Strang-Fix conditions.

\begin{defi}\label{d2}
We say that a  function $\phi$ satisfies {\em the Strang-Fix conditions} of
order $s$  if $D^{\beta}\h{\phi}(k) = 0$ for every $\beta\in\z_+^d$, $[\beta]<s$,
and for all $k\in\zd\setminus \{\nul\}$.
\end{defi}

Certain compatibility conditions for a distribution $\w\phi$ and a function $\phi$
is also required to provide good approximation properties of the operator $Q_j(f,\phi,\w\phi)$. For our purposes, we will use the following two conditions.

\begin{defi}
\label{d3}
 A tempered distribution  $\w\phi$ and a function $\phi$ is said to be  {\em weakly compatible of order~$s$} if $D^{\beta}(1-\h\phi\h{\w\phi})({\bf 0}) = 0$ 	for every $\beta\in\z_+^d$, $[\beta]<s$.
\end{defi}

\begin{defi}
\label{d1}
 A tempered distribution  $\w\phi$ and a function $\phi$ is said to be  {\em strictly compatible} if there exists $\delta>0$ such that  $\overline{\h\phi}(\xi)\h{\w\phi}(\xi)=1$  a.e. on $\delta\td$.
\end{defi}

For $\w\phi\in \mathcal{S}'_N$ and different classes
of  functions $\phi$, approximation properties of quasi-projection operators
$
Q_j(f,\phi,\w\phi)%=\sum_{k\in\zd} \langle f, {\w\phi}_{jk}\rangle \phi_{jk},
$
were studied  in~\cite{Sk1},  \cite{KS}, \cite{KS20}, and~\cite{KKS}. Generally speaking, if $\w\phi\in \mathcal{S}'$, then the functional $\langle f, \w\phi_{jk}\rangle $ has meaning only for functions $f$ belonging $\mathcal{S}$.
Under some additional restrictions on the distribution $\w\vp$, the class of functions $f$ can be essentially extended. To this end, the quantity $\langle f, {\w\phi}_{jk}\rangle $ was replaced  by the inner product $\langle \h f, \h{\w\phi_{jk}}\rangle$ in the mentioned papers.
The following result  is a combination of Theorem~14 from~\cite{cksv} and Theorem~5 from~\cite{KS}.

\medskip

\noindent{\bf Theorem A.}	
%\label{theoQj}
	{\it Let $2\le p < \infty$, $s\in\n$,  $N\ge0$,  $\delta\in(0, 1/2)$,  $M$ be an isotropic matrix, $\psi \in L_p$, and $\w\psi \in \mathcal{S}_N'$.
	Suppose
\begin{itemize}
 	\item[$1)$]  $\h\psi\in L_{p'}$ and
	%there exists $B_{\psi}>0$ such that
$\sum\limits_{k\in\zd}  |\h\psi(\xi+k)|^{p'}<c_1$ for all $\xi\in\rd$;

      \item[$2)$]  $\h\psi(\cdot+l)\in C^{s}(B_\delta)$  for all $l\in\zd\setminus\{\nul\}$
and $\sum\limits_{l\ne\nul}\sum\limits_{\|\beta\|_1=s}
		\sup\limits_{|\xi|<\delta}|D^\beta\h\psi(\xi+l)|^{p'}< c_2$;

	 \item[$3)$]  the Strang-Fix conditions of order $s$ are satisfied	for $\psi$;
	
      \item[$4)$]  $\overline{\h\psi}\h{\w\psi}\in C^s(B_\delta)$;
 	
   \item[$5)$]  $\psi$ and ${\w\psi}$ are weakly  compatible of order~$s$.
\end{itemize}

\noindent
 		If $\h f\in L_{p'}$, and
 $\h f(\xi)=\mathcal{O}(|\xi|^{-N-d-\varepsilon})$, $\varepsilon>0$,
as $|\xi|\to\infty$, then

	 $$
\bigg\|f-\lim\limits_{N\to \infty}\sum_{\| k\|_\infty \le N} \langle \widehat{f}, \h{\w\psi_{jk}}\rangle \psi_{jk}\bigg\|_p\le C
	 \begin{cases}
	 |\lambda|^{-j(N+d/p + \varepsilon)}  &\mbox{if }
	s> N+d/p+ \varepsilon\\
	  (j+1)^{1/p'} |\lambda|^{-js} &\mbox{if }
	 s= N+d/p + \varepsilon \\
	|\lambda|^{-js}
	 &\mbox{if }
	 s< N+d/p + \varepsilon
	\end{cases},
$$
	where $\lambda$ is an eigenvalue of $M$ and the constant $C$ is independent on $j$.
}
%\end{theoA}	

\medskip

This result is obtained for a wide class of  operators $Q_j(f,\psi,\w\psi)$, but
unfortunately, the error estimate is given only for $p\ge2$. Another drawback of this theorem
is the restriction on  the decay of  $\h f$. It is not difficult to see that it is redundant, for example,
if $\w\psi \in L_{p'}$ and $f\in L_p$. Also, although  Theorem~A provides approximation order for $Q_j(f,\psi,\w\psi)$,  more accurate  error estimates in terms of smoothness of $f$ were not obtained.

The mentioned drawbacks of Theorem~A were avoided
in~\cite{KS3}, where a class of Kantorovich-type  operators $Q_j(f,\phi,\w\phi)$ associated with a regular function $\w\phi$
and a bandlimited function $\phi$ was considered. In particular,  the next theorem was obtained in~\cite[Theorem~17]{KS3}. To formulate it, we introduce the space ${\cal B}$, which consists
of  functions $\phi$ given by
$
\phi=\F^{-1}\theta,
$
  where $\supp\theta\subset [a,b]:=[a_1, b_1]\times\dots\times[a_d, b_d]$ and
	$\theta\big|_{[a,b]}\in C^d([a,b])$.

\medskip

\noindent{\bf Theorem B.}
%\begin{theo}
%\label{KSth2}
{\it Let $1<p<\infty$,  $s\in\n$, $\delta>0$, and $\varepsilon\in (0,1)$.  Suppose
\begin{itemize}
  \item[$1)$] $\phi\in\cal B$, $\supp\h\phi\subset B_{1-\varepsilon}$, and $\h\phi\in C^{s+d+1}(B_\delta)$;
%either $\w\phi\in L_q^0$,  or
 \item[$2)$]  $\w\phi \in\mathcal{B}\cup \mathcal{L}_{p'}$ and  $\h{\w\phi}\in C^{s+d+1}(B_\delta)$;
 \item[$3)$] $\phi$ and ${\w\phi}$ are weakly  compatible of order~$s$.
\end{itemize}
 Then, for every $f\in L_p$, we have
\begin{equation*}
  \bigg\|f-\sum\limits_{k\in\zd}
\langle f,\widetilde\phi_{jk}\rangle \phi_{jk}\bigg\|_p\le c\,\omega_s\(f,\|M^{-j}\|\)_p,
\end{equation*}
where $c$ is independent on $f$ and $j$.
}
%\end{theo}

\medskip

In what follows, we will consider  a class of quasi-projection operators $Q_j(f,\vp,\w\vp)$ associated with a tempered distribution $\w\phi$ belonging to the class $\mathcal{S}_{\a,p;M}'$, where  $1\le p\le \infty$, $M\in \mathfrak{M}$, and $\a\in \mathcal{A}_M$. We say that  $\w\phi\in\mathcal{S}_{\a,p;M}'$  if $\h{\w\phi}$ is a measurable locally bounded function  and
\begin{equation}\label{DefS}
  \Vert \Lambda_{\F({\w\phi^-})}P_\mu \Vert_p \le \alpha (M^\mu) \Vert P_\mu \Vert_p\quad\text{for all}\quad P_\mu\in \mathcal{B}_{M^\mu,p}\cap L_2,\,\, \mu\in \Z_+.
\end{equation}
%for all $\mu\in \Z_+$ and $T_\mu\in \mathcal{B}_{M^\mu,p}\cap L_2$.

Obviously, inequality~\eqref{DefS} is satisfied  with $\alpha\equiv 1$  if $\w\phi$ is the Dirac delta-function or $\w\vp\in L_1$.
If $M={\rm diag}(m_1,\dots,m_d)$ and $\w\phi$ is a distribution corresponding to
the differential operator of the form $\w\vp(x)=D^\b \d(x)$, $\beta\in \z_+^d$,
then $\w\phi$ belongs to the class $\mathcal{S}_{\a,p;M}'$ with $\alpha(M)=m_1^{\beta_1}\dots m_d^{\beta_d}$.
If $M$ is an isotropic matrix, then $\alpha(M)=m^{{[\beta]}/d}$.
This follows from the Bernstein inequality (see, e.g.,~\cite[p.~252]{Timan}) given by
$$
\Vert P'\Vert_{L_p(\R)} \le \sigma \Vert P \Vert_{L_p(\R)}, \quad P\in L_p(\R),\quad \supp \h P \subset [-\sigma,\sigma].
$$

Now we are going to extend the operator $Q_j(f,\phi,\w\phi)$ with $\w\phi \in \mathcal{S}_{\a,p;M}'$
onto the Besov spaces $\mathbb{B}_{p;M}^{\a(\cdot)}$ and the space ${C}_0$. For this, we need to
define (extend) the functional $\langle f,\w\phi_{jk}\rangle$ in an appropriate way. In the case
$\a(M)=|\det M|^{N/d+1/p}$ and $1\le p<\infty$, a similar extension was given in~\cite{KS19}.
	
\begin{defi}
\label{def0}
	Let $1\le p\le\infty$,  $M\in \mathfrak{M}$, $\a\in \mathcal{A}_M$, and $\d\in(0,1]$. For
	$\w\vp\in \mathcal{S}_{\a,p;M}'$
and $f\in \mathbb{B}_{p;M}^{\a(\cdot)}$  or $\w\vp\in \mathcal{S}_{{\rm const},\infty;M}'$ and $f\in {C}_0$,   we set
\begin{equation}\label{con}
  \langle f,\w\vp_{0k}\rangle:%=\lim_{\mu\to\infty} P_\mu*[\mathcal{F}^{-1}(\eta({M^{*-\mu-\mu_0}}\cdot))*\w\vp^-](k)
  =\lim_{\mu\to\infty} \langle \h{P_\mu},\h{\w\vp_{0k}}\rangle, \quad k\in\zd,
\end{equation}
where the functions $\{P_\mu\}_{\mu\in\Z_+}$ are such that $P_\mu\in \mathcal{B}_{\d M^\mu,p}\cap L_2$  and
\begin{equation}\label{eT}
  \|f-P_\mu\|_p\le c(d,p){E}_{\d_p M^\mu}(f)_p, \quad \delta_p=\left\{
                                                 \begin{array}{ll}
                                                \delta & \hbox{if $ p<\infty$,} \\
                                                   \delta/2 & \hbox{if $p=\infty$.}
                                                 \end{array}
                                               \right.
\end{equation}
Set also
$$
\langle f,\w\vp_{jk}\rangle:= m^{-j/2}\langle f(M^{-j}\cdot),\w\vp_{0k}\rangle, \quad j\in\z_+.
$$
\end{defi}

Some comments are needed to approve this definition.
First, it will be proved in  Lemma~\ref{lem1} that   the limit in~\eqref{con} exists and does not depend on a choice of  $P_\mu$ and  $\delta$.
Second, in view of Lemmas~\ref{lem2} and~\ref{lemJ2}, one can always find  functions $P_\mu\in \mathcal{B}_{\d M^\mu,p}\cap L_2$, $1\le p\le\infty$,  such that~\eqref{eT} holds. Third, we can  write
\begin{equation}\label{conv}
  \langle \h{P_\mu},\h{\w\vp_{0k}}\rangle=\Lambda_{\mathcal{F}(\w\phi^-)}P_\mu(-k).%=T_\mu*(\F^{-1}(\eta(M^{*-\mu}\cdot))*\w\vp^-)(-k).
\end{equation}
Finally, we mention that if $\w\vp \in L_{p'}$, then $\langle f,\w\phi_{jk}\rangle$ is the standard inner product, which has meaning for any $f\in L_p$.
	
\begin{rem}
\label{prop003}
%If the Fourier transform of $f$ has enough
%decay for the inner product $\langle \h f,\h{{\w \vp}_{0k}}\rangle$ to have meaning, then
%it is natural to define $Q_j(f,\phi,\w\phi)$ setting $\langle f,\w\vp_{0k}\rangle:=\langle \h f,\widehat{{\w \vp}_{0k}}\rangle$
%(see Theorem~A).  Such an operator $Q_j(f,\phi,\w\phi)$  will be the same as the one in the correspondence of
%Definition~\ref{def0}. The explanation of this fact can be found in~\cite{KS19}.
%In the series of works~\cite{??}, we studied properties of the operator $Q_j(f,\phi,\w\phi)$ for distributions $\w\vp\in \mathcal{S}_N'$ and functions $f$ with a good decay of its Fourier transform (see also Theorem~A).

Note that if $\w\vp\in \mathcal{S}_N'$ for some $N\ge 0$ and the Fourier transform of a function $f$ has a sufficiently good decay such that the inner product $\langle \h f,\h{{\w \vp}_{0k}}\rangle$ has sense, then it is natural to define the operator $Q_j(f,\phi,\w\phi)$ by setting $\langle f,\w\vp_{0k}\rangle:=\langle \h f,\widehat{{\w \vp}_{0k}}\rangle$ (see, e.g.,~\cite{KS}, \cite{KS20}, \cite{KKS} as well as Theorem~A). It is not difficult to see that such an operator $Q_j(f,\phi,\w\phi)$ is the same as the corresponding operator defined by means of
Definition~\ref{def0} (see, e.g.,~\cite{KS19}). %The explanation of this fact can be found in~\cite{KS19}.
\end{rem}

The main tools in this paper are Fourier multipliers. Let us recall their definition and basic properties.
\begin{defi}\label{dmult}
Let $h$ be a bounded measurable function on $\R^d$. Consider the linear transformation  $\Lambda_h$ defined by
$\Lambda_h (f)= \mathcal{F}^{-1}(h \h f )$, $f\in L_2\cap L_p$. The function $h \,:\, \R^d \to \mathbb{C}$ is called a Fourier multiplier in $L_p$, $1\le p\le \infty$, (we write $h \in \mathcal{M}_p$) if there exists a constant $K$ such that
\begin{equation}\label{nmul}
  \Vert \Lambda_h (f)\Vert_p\le K \Vert f\Vert_p\quad\text{for any}\quad f\in L_2\cap L_p.
\end{equation}
The smallest $K$, for which inequality~\eqref{nmul} holds, is called the norm of the multiplier $h$. We denote this  norm by  $\Vert h \Vert_{\mathcal{M}_p}$.
\end{defi}

Note that if~\eqref{nmul} holds and $1\le p<\infty$, then the operator $\Lambda_h$ has a unique bounded extension to $L_p$, which satisfies the same inequality. As usual, we denote this extension by~$\Lambda_h$.

Let us recall some basic properties of Fourier multipliers
(see, e.g.,~\cite[Ch.~6]{BL} and~\cite[Ch.~1]{Nik}):

\begin{itemize}
  \item[$(i)$] if $1<p<2$, then $\mathcal{M}_1\subset
\mathcal{M}_p\subset \mathcal{M}_2=L_\infty$;

\item[$(ii)$] if $1\le p\le\infty$, then
$\mathcal{M}_p=\mathcal{M}_{p'}$ and $\Vert h\Vert_{\mathcal{M}_{p}}=\Vert h\Vert_{\mathcal{M}_{p'}}$;

\item[$(iii)$] if $h_1, h_2\in \mathcal{M}_p$, then $h_1+h_2\in
\mathcal{M}_p$ and $h_1h_2\in \mathcal{M}_p$;

\item[$(iv)$] if $h\in \mathcal{M}_p$,  then
$h(A\cdot)\in \mathcal{M}_p$ and $\|h(A\cdot)\|_{\mathcal{M}_p}=\|h \|_{\mathcal{M}_p}$ for any non-singular matrix $A$.
\end{itemize}

%\medskip

The classical sufficient condition for Fourier multipliers in $L_p$, $1<p<\infty$, is Mikhlin's condition (see, e.g.,~\cite[p.~367]{Gr}), which states that
{if a function $h$ is such that
\begin{equation*}\label{VV.eq1}
    |D^\nu h(\xi)|\le K|\xi|^{-[\nu]},\quad \xi\in\R^d\setminus \{0\},
\end{equation*}
for all $\nu\in \Z_+^d$, $[\nu]\le d/2+1$, then
$h\in \mathcal{M}_p$ for all $1<p<\infty$ and $\Vert h\Vert_{\mathcal{M}_p}\le c(p,d)\(\Vert h\Vert_\infty+K\)$.}

Concerning the limiting cases $p=1$ and $\infty$, we  note that if $h$ is a continuous function, then $h\in \mathcal{M}_1$ if and only if $h$ is the Fourier
transform of a finite Borel (complex-valued) measure. The multiplier itself is a convolution of a function and this measure. Numerous efficient sufficient conditions for Fourier multipliers in $L_1$ and $L_\infty$ can be found in the survey~\cite{LST}. Here, we only mention the Beurling-type condition, which states that if $h\in W_2^{k}$ with $k>d/2$, then $h\in \mathcal{M}_1$ (see, e.g.,~\cite[Theorem~6.1]{LST}).

Finally, we note that if  $\vp\in \mathcal{B}$, then $\h\vp\in \mathcal{M}_p$ for all $1<p<\infty$. Indeed, for any $\vp\in \mathcal{B}$,  we have $\h\vp=\chi_\Pi \cdot \theta$, where $\theta$ belongs to $C^d(\R^d)$  and has a compact support. It is well known that the characteristic function of $\Pi$ is a Fourier multiplier in $L_p$, $1<p<\infty$ (see, e.g.,~\cite[p.~100]{Stein}). By Mikhlin's condition the same holds for the function $\theta$. Thus, it follows from $(iii)$ that $\h\vp\in \mathcal{M}_p$.

\section{Auxiliary results}
	
\begin{lem}\label{lemMZ} {\sc (\cite[Theorem 4.3.1]{TB})}
Let $g\in L_p$, $1\le p<\infty$, and $\supp \h g\subset [-\sigma_1,\sigma_1]\times\dots\times [-\sigma_d,\sigma_d]$, $\sigma_j>0$, $j=1,\dots,d$. Then
$$
\frac1{\sigma_1\dots \sigma_d} \sum_{k\in \Z^d} \max_{x\in Q_{k,\sigma}} |g(x)|^p \le c\,\Vert g\Vert_p^p,
$$
where $Q_{k,\sigma}=[\frac{2k_1-1}{2\sigma_1}, \frac{2k_1+1}{2\sigma_1}]\times\dots\times [\frac{2k_d-1}{2\sigma_d}, \frac{2k_d+1}{2\sigma_d}]$ and $c$ depends only on $p$ and $d$.  %does not depend on $j$.
\end{lem}

\begin{lem}\label{lem001}
Let  $1\le p\le\infty$, $g\in L_p$,   $h\in L_{p'}$, and  $\h h\in \mathcal{M}_{p}$. Then the operator  $T(g):=h*g$ is bounded in $L_p$ and $\|h*g\|_p\le \|\h h\|_{\mathcal{M}_p}\|g\|_p$.
\end{lem}
{\bf Proof.} In the case $p=\infty$, the statement follows from Minkowski's inequality (without assumption $\h h\in \mathcal{M}_{p}$).
 Consider the case $p<\infty$. Choose a sequence $\{g_n\}_n\subset \mathcal{S}$ converging to $g$ in $L_p$-norm. Since $\h h\in \mathcal{M}_{p}$,  the functions $\Lambda_{\h h}(g_n)$ form a Cauchy sequence  in $L_p$.
Hence, $\Lambda_{\h h}(g_n)\to G$, $G \in L_p$.
On the other hand, $\Lambda_{\h h}(g_n)=h*g_n$, and the sequence $h*g_n$ converges to $h*g$ almost everywhere. It follows that
 $(h*g)(x)=G(x)$ for almost all $x$. Thus, we derive
\begin{equation*}
  \begin{split}
      \|h*g\|_p&=\|G\|_p\le \|\Lambda_{\h h}g_n\|_p+\|G-\Lambda_{\h h}g_n\|_p\le \|\h h\|_{\mathcal{M}_p}\|g_n\|_p+\|G-\Lambda_{\h h}g_n\|_p\\
      &\le \|\h h\|_{\mathcal{M}_p}\|g\|_p+\|\h h\|_{\mathcal{M}_p}\|g-g_n\|_p+\|G-\Lambda_{\h h}g_n\|_p.
   \end{split}
\end{equation*}
Finally, passing to the limit as $n\to\infty$, we complete the proof.~~$\Diamond$

\begin{lem}
\label{lem1}
Let  $1\le p\le \infty$, $M\in \mathfrak{M}$, $n\in\N$, $\delta\in (0,1]$, and $\a\in \A_M$.
Suppose that $\w\vp$, $f$, and the functions $P_\mu$, $\mu\in\Z_+$, are as in Definition~\ref{def0}.
Then the sequence $\{\{\langle \h{P_\mu},\h{\w\phi_{0k}}\rangle \}_k\}_{\mu=1}^\infty$ converges in $\ell_p$ as $\mu\to\infty$ and its limit is  independent of the choice of $P_\mu $
and $\delta$; a fortiori for every $k\in\zd$ there exists the limit
$\lim_{\mu\to\infty}\langle \h{P_\mu},\h{\w\phi_{0k}}\rangle $  independent of the choice of $P_\mu $ and $\delta$.  Moreover,
for all $f\in \mathbb{B}_{p;M}^{\a(\cdot)}$,
we have
\begin{equation}\label{elem1}
	\sum_{\mu=n}^\infty
	\|\{\langle \h{P_{\mu+1}},\h{\w\phi_{0k}}\rangle -\langle \h{P_\mu},\h{\w\phi_{0k}}\rangle \}_{k}\|_{\ell_p}\le
	 c\sum_{\mu=n}^\infty m^\frac\mu p \alpha (M^\mu)
E_{\d_p M^\mu}(f)_p,
\end{equation}
where $c$ depends only on $d$, $p$, and $M$.
\end{lem}
{\bf Proof.}
Consider the case $p<\infty$.
Setting
$$
F(x):=\int_{\rd}\Big(\h{P_{\mu+1}}({M^*}^{\mu+1}\xi)-\h{P_{\mu}}({M^*}^{\mu+1}\xi)\Big)
\overline{\h{\w \phi}({M^*}^{\mu+1}\xi)}e^{2\pi i(\xi,x)}\,d\xi,
$$
we get
\begin{equation}
\label{121}
\begin{split}
     &\|\{\langle\h{P_{\mu+1}}, \h{\w\phi_{0k}}\rangle -\langle \h{P_\mu},\h{\w\phi_{0k}}\rangle \}_{k}\|_{\ell_p}^p
     =m^{p(\mu+1)}\sum_{k\in\zd}|F(M^{\mu+1}k)|^p.
\end{split}
\end{equation}
Since
$\supp \h F\subset [-\sigma, \sigma]^d$, where  $\sigma=\sigma(M,d)>1$,
using Lemma~\ref{lemMZ} and  taking into account that each set $Q_{k,\sigma}$
contains a finite number (depending only on $M$ and $d$) points $M^{\mu+1}k$,  $k\in\zd$, we obtain
	\be
\label{122}
	\sum_{k\in\zd}|F(M^{\mu+1}k)|^p\le c_1\sigma^d\int_{\rd}|F(x)|^p\,dx.
	\ee
Recall that $\mu_0=\min\{\mu\in \N\,:\,\T^d\subset \frac12 M^{*\nu}\T^d\quad\text{for all}\quad\nu\ge \mu-1\}$. Since
${M^*}^{\mu+1}\delta\td
\subset {M^*}^{\mu+\mu_0}\delta\td$ and ${M^*}^{\mu}\delta\td
\subset {M^*}^{\mu+\mu_0}\delta\td$, both the functions
$P_\mu$ and $P_{\mu+1}$ are in $\mathcal{B}_{\d M^{\mu+\mu_0},p}\cap L_2$.
Thus, combing~\eqref{121} and~\eqref{122} and using~\eqref{DefS},  we derive
\begin{equation}\label{121+}
  \begin{split}
    \|\{\langle\h{P_{\mu+1}}, \h{\w\phi_{0k}}\rangle -&\langle \h{P_\mu},\h{\w\phi_{0k}}\rangle \}_{k}\|_{\ell_p}\le c_2m^{\mu+1}\|F\|_p=c_2m^{\mu+1+\frac{\mu+1}{p}}\|F(M^{\mu+1}\cdot)\|_p\\
    &=c_2m^{\frac{\mu+1}{p}}\|\L_{\F(\w\vp^-)}(P_{\mu+1}-P_{\mu})\|_p\le c_2 m^{\frac{\mu+1}{p}}\a(M^{\mu+\mu_0})\|P_{\mu+1}-P_{\mu}\|_p\\
    &\le c_3 m^{\frac \mu p} \big(\alpha (M^{\mu}) E_{\delta M^\mu}(f)_p + \alpha (M^{\mu+1}) E_{\delta M^{\mu+1}}(f)_p\big),
  \end{split}
\end{equation}
which implies~\eqref{elem1} after the corresponding summation.

Let now $p= \infty$. Taking into account~\eqref{conv},
we can write
$$
\langle\h{P_{\mu+1}}, \h{\w\phi_{0k}}\rangle -\langle \h{P_\mu},\h{\w\phi_{0k}}\rangle=\Lambda_{\mathcal{F}(\w\phi^-)}(P_{\mu+1}-P_\mu)(-k).
$$
Then, using~\eqref{DefS}, we obtain
\begin{equation}\label{121++}
  \begin{split}
      \|\langle\h{P_{\mu+1}}, \h{\w\phi_{0k}}\rangle -\langle \h{P_\mu},\h{\w\phi_{0k}}\rangle\|_{\ell_\infty}&\le\|\Lambda_{\mathcal{F}(\w\phi^-)}(P_{\mu+1}-P_\mu)\|_\infty\le \a(M^{\mu+\mu_0})\|P_{\mu+1}-P_\mu\|_\infty\\
      &\le c_4\big(\alpha(M^{\mu}) E_{\frac{\d}2 M^\mu}(f)_\infty +  \alpha (M^{\mu+1})E_{\frac{\d}2 M^{\mu+1}}(f)_\infty\big),
   \end{split}
\end{equation}
which again implies~\eqref{elem1}.

Next, it is clear that there exists $\nu(\d)\in \N$ such that $E_{\delta_p M^\mu}(f)_p\le E_{M^{\mu-\nu(\delta)}}(f)_p$ and $\alpha(M^{\mu})\le c(\delta,M)\alpha(M^{\mu-\nu(\delta)})$
 for all big enough $\mu$.
Thus, if  $f\in \mathbb{B}_{p;M}^{\a(\cdot)}$, then it follows from ~\eqref{elem1} that
$\{\{\langle \h{P_\mu},\h{\w\phi_{0k}}\rangle\}_k\}_{\mu=1}^\infty$ is a Cauchy sequence in
$\ell_p$.  Fortiori, for every $k\in\zd$, the sequence $\{\langle \h{P_\mu},\h{\w\phi_{0k}}\rangle\}_{\mu=1}^\infty$  has a limit.

Let now $p=\infty$, $\alpha=\const$, and $f\in {C}_0$. For every $\mu', \mu''\in\n$, there exists
$\nu\in\n$ such that both the functions $\h{P_{\mu'}}$ and $\h{P_{\mu''}}$ are supported
in $M^{*\nu}\T^d$, and similarly to~\eqref{121++}, we have
$$
 \|\{\langle\h{P_{\mu'}}, \h{\w\phi_{0k}}\rangle -\langle \h{P_{\mu''}},\h{\w\phi_{0k}}\rangle\}_{k}\|_{\ell_\infty}\le c_5\big(E_{\frac{\d}2 M^{\mu'}}(f)_\infty +
 E_{\frac{\d}2 M^{\mu''}}(f)_\infty\big).
$$
Thus, again $\{\{\langle\h{P_{\mu}}, \h{\w\phi_{0k}}\rangle\}_k\}_{\mu=1}^\infty$ is a Cauchy sequence in
$\ell_\infty$ and  every  sequence $\{\langle\h{P_{\mu}}, \h{\w\phi_{0k}}\rangle\}_{\mu=1}^\infty$  has a limit.

Let us  check that the limit of $\{\{\langle\h{P_{\mu}}, \h{\w\phi_{0k}}\rangle\}_k\}_{\mu=1}^\infty$ in $\ell_p$ does not depend on the choice of functions $P_\mu$ and $\delta$.  Let $\delta'\in (0,1]$ and
$P'_\mu\in \mathcal{B}_{\d' M^\mu,p}\cap L_2$  be such that
$  \|f-P'_\mu\|_p\le c'(d,p){E}_{\d'_p M^\mu}(f)_p$.
Since both the functions
$P_\mu$ and $P'_{\mu}$ are in $\mathcal{B}_{M^{\mu},p}\cap L_2$, repeating the
arguments of the proof of inequalities~\eqref{121+} and~\eqref{121++} with $P'_\mu$ instead of $P_{\mu+1}$ and
$0$ instead of $\mu_0$, we obtain
\begin{equation*}
  \begin{split}     	
	\|\{\langle\h{P_{\mu}'}, \h{\w\phi_{0k}}\rangle -\langle \h{P_{\mu}},\h{\w\phi_{0k}}\rangle\}_{k} \|_{\ell_p}
&\le c_6\,m^{\frac \mu p}\alpha(M^\mu) \|P'_{\mu}-P_\mu\|_p\\
   &\le  c_7\,m^{\frac \mu p} \alpha (M^{\mu}) (E_{\delta_p M^\mu}(f)_p+E_{\delta'_p M^\mu}(f)_p ).
   \end{split}
\end{equation*}
It follows that
$\|\{\langle\h{P_{\mu}'}, \h{\w\phi_{0k}}\rangle -\langle \h{P_{\mu}},\h{\w\phi_{0k}}\rangle\}_{k} \|_{\ell_p}\to 0$
as $\mu\to\infty$, which yields the independence from the choice of $P_\mu $ and $\delta$.~$\Diamond$

\begin{lem}\label{lemKK1+}
Let  $\phi\in L_p$, $1\le p\le\infty$, be such that $\supp \h\vp$ is compact and
$\h\phi\in \mathcal{M}_{p}$. Then, for any sequence $\{a_k\}_{k\in\zd}\in \ell_p$ if $p<\infty$ and $\{a_k\}_{k\in\zd}\in {\rm c_0}$ if $p=\infty$, the series $\sum_{k\in \Z^d} a_k \phi_{0k}(x)$ converges  unconditionally in $L_p$ and
\begin{equation*}
%\label{lemsum}
  \bigg\Vert \sum_{k\in \Z^d} a_k \phi_{0k}\bigg\Vert_p
\le  c\left\Vert \{a_k\}_{k}\right\Vert_{\ell_{p}},
\end{equation*}
where  $c$ does not depend on $\{a_k\}_{k\in \Z^d}$.
\end{lem}
{\bf Proof.}
Let us fix an integer $n$.
 By duality, we can find a function  $g\in L_{p'}$  such that $\Vert g\Vert_{p'}\le 1$ and
\begin{equation}\label{eqKK1}
\begin{split}
   \bigg\Vert \sum_{\|k\|_\infty\le n} a_k \phi(\cdot-k)\bigg\Vert_p&=\bigg|\bigg\langle \sum_{\|k\|_\infty\le n} a_k \phi(\cdot-k),g\bigg\rangle\bigg|=\bigg|\sum_{\|k\|_\infty\le n} a_k \langle \phi(\cdot-k),g \rangle\bigg|.
\end{split}
\end{equation}

 Consider the case $p>1$. Applying H\"older's inequality, using Lemmas~\ref{lemMZ} and \ref{lem001}, and  taking into account that
$\mathcal{M}_{p}=\mathcal{M}_{p'}$, we obtain
\begin{equation}\label{eqKK2}
\begin{split}
\sum_{\|k\|_\infty\le n}  &|a_k \langle \phi(\cdot-k),g \rangle|\le\left\Vert \{a_k\}_{k}\right\Vert_{\ell_{p}}
\bigg(\sum_{k\in \Z^d} |\langle \phi(\cdot-k),g \rangle|^{p'} \bigg)^\frac1{p'}\\
&=\left\Vert \{a_k\}_{k}\right\Vert_{\ell_{p}}\bigg(\sum_{k\in \Z^d} |(\phi*g^-)(k)|^{p'} \bigg)^\frac1{p'}\le c_1\left\Vert \{a_k\}_{k}\right\Vert_{\ell_{p}}\Vert\phi*g^-\Vert_{p'}\le c_2\left\Vert \{a_k\}_{k}\right\Vert_{\ell_{p}}
\Vert g \Vert_{p'},
\end{split}
\end{equation}
where $c_2$ does not depend on $n$.
Combining~\eqref{eqKK1} and~\eqref{eqKK2}, we get  that the cubic sums of the series $\sum_{k\in \Z^d} a_k \phi_{0k}$
are bounded in  $L_p$-norm.

Similarly,  the  boundedness of  the cubic sums in  $L_1$-norm follows from
\begin{equation*}
  \begin{split}
      \sum_{k\in \Z^d} |a_k \langle \phi_j(\cdot-k),g \rangle|&\le\left\Vert \{a_k\}_{k}\right\Vert_{\ell_{1}}
\sup_{k}|\langle \phi(\cdot-k),g \rangle|\\
&\le c_3\left\Vert \{a_k\}_{k}\right\Vert_{\ell_{1}}\Vert\phi*g^-\Vert_\infty\le c_4\left\Vert \{a_k\}_{k}\right\Vert_{\ell_{1}}
\Vert g \Vert_\infty.
   \end{split}
\end{equation*}
Now it is clear that all statements hold.~~$\Diamond$

\begin{lem} {\sc (\cite[Theorem 2.1]{v58})}
\label{lemKK1}
Let $\phi\in {\cal L}_p$, $1\le p\le\infty$. Then, for any sequance $\{a_k\}_{k\in\zd}\in\ell_p$, we have
\begin{equation*}
%\label{lemsum}
  \bigg\Vert \sum_{k\in \Z^d} a_k \phi_{0k}\bigg\Vert_p
\le  \|\phi\|_{\mathcal{L}_p} \left\Vert \{a_k\}_{k}\right\Vert_{\ell_{p}}.
\end{equation*}
\end{lem}

\begin{lem} {\sc (\cite[Proposition 5]{NU})}
\label{lemNU}
Let $f\in L_p$,  $1\le p\le\infty$, and  $\w\phi\in {\cal L}_{p'}$. Then
\begin{equation*}
%\label{lemsum}
\left\Vert \{\langle f, \w\phi_{0k}\rangle\}_{k}\right\Vert_{\ell_{p}}
\le  \|\w\phi\|_{\mathcal{L}_{p'}}\|f\|_p .
\end{equation*}
\end{lem}

\begin{lem}
\label{lem99}
Let $f\in {C}_0$ and $\w\vp\in \mathcal{S}_{\const,\infty;M}'$ for some $M\in \mathfrak{M}$. Then $\{\langle f, \w\phi_{0k}\rangle\}_k\in {\rm c}_0$ and
$$
\|\{\langle f, \w\phi_{0k}\rangle\}_k\|_{\ell_\infty}\le c\,\|f\|_\infty,
$$
where $c$ does not depend on $f$.
\end{lem}
{\bf Proof.} By Lemma~\ref{lem1}, for any $\varepsilon>0$, there exists a function   $P_\mu\in \mathcal{B}_{M^\mu,\infty}\cap L_2$
such that   $\|f-P_\mu\|_\infty\le c_1{E}_{M^\mu}(f)_\infty$ and
\begin{equation}\label{yu1}
  \|\{\langle f, \w\phi_{0k}\rangle-\langle \h{P_\mu}, \h{\w\phi_{0k}}\rangle\}_k\|_{\ell_\infty}<\varepsilon.
\end{equation}
Since $\h{P_\mu}\in L_2$ has
a compact support and $\h{\w \phi}$ is locally bounded, the function $\h{P_\mu}\h{\w \phi}$ is summable. Hence,
$
\langle \h{P_\mu}, \h{\w\phi_{0k}}\rangle=
\mathcal{F}(\h{P_{\mu}}\,\h{\w \phi^-})(k)\to 0
$
as $|k|\to \infty$. Thus, using~\eqref{yu1}, we get that $\langle f, \w\phi_{0k}\rangle\to 0$ as $|k|\to \infty$. Moreover, due to~\eqref{DefS} and~\eqref{conv}, we have
$$
|\langle \h{P_\mu}, \h{\w\phi_{0k}}\rangle|=
|\Lambda_{\F({\w\phi^-})}P_{\mu} (-k)|
\le\Vert \Lambda_{\F({\w\phi^-})}P_{\mu} \Vert_\infty \le \const\Vert P_{\mu}\Vert_\infty
\le c_1\|f\|_\infty.
$$
Combining this with~\eqref{yu1}, we prove the lemma.~$\Diamond$

\medskip

Let us also recall some basic inequalities for the best approximation and moduli of smoothness.

\begin{lem}{\sc (\cite[Lemma 8]{KS19})}\label{lem2}
Let $f\in L_p$,  $1\le p< \infty$, and let $A$ be a $d\times d$-matrix.
Then
$$
\inf_{P\in \mathcal{B}_{A,p}\cap L_2}\|f-P\|_p\le c\,E_{A} (f)_p,
$$
where $c$ depends only on $p$ and $d$.
\end{lem}
	
\begin{lem}\label{lemJ2}
Let $f\in {C}_0$ and let $A$ be a $d\times d$-matrix. Then
%\be
%\label{201}
$$
\inf_{P\in \mathcal{B}_{2A,\infty}\cap L_2}\|f-P\|_\infty\le c\,E_{A} (f)_\infty,%\quad E_{I}(f)_\infty\le C\omega_s(f,1)_\infty,
$$
%\ee
where $c$ depends only on $d$.
\end{lem}
{\bf Proof.}
Since $f\in {C}_0$, there exists a compactly supported $g\in {C}$ satisfying
$$
\|f-g\|_\infty\le E_{A} (f)_\infty.
$$
Let $Q\in \mathcal{B}_{A,\infty}$ be such that
$$
\|g-Q\|_\infty\le 2E_{A} (g)_\infty.
$$
Denote $N_A=\mathcal{F}^{-1}(\eta(A^{*-1}\cdot))$.   Obviously, $N_A*Q=Q$ and
$N_A*g\in \mathcal{B}_{2A,\infty}\cap L_2$. This, together with the above two inequalities, yields
\begin{equation*}
  \begin{split}
      \inf_{P\in \mathcal{B}_{2A,\infty}\cap L_2}\|f-P\|_\infty&\le\|f-N_A*g\|_\infty\le
\|f-g\|_\infty+\|g-Q\|_\infty+\big\| N_A*(g-Q)\big\|_\infty\\
&\le 3E_{A} (f)_\infty+\|N_A\|_1\|g-Q\|_\infty\le
(3+\|N_I\|_1) E_{A} (f)_\infty.~~\Diamond
   \end{split}
\end{equation*}

%\medskip

\begin{lem}\label{lemJ} {\sc (See~\cite[5.2.1 (7)]{Nik} or~\cite[5.3.3]{Timan})}
Let $f\in L_p$, $1\le p \le\infty$, and $s\in \N$. Then
$$
E_{I}(f)_p\le c\,\omega_{s}(f, 1)_p,
$$
where $c$ depends only on $d$ and $s$.
\end{lem}

Finally, the next two statements can be found, e.g., in ~\cite{Wil}, see also~\cite{KT20}.
\begin{lem}\label{lemNS}
  Let $P\in \mathcal{B}_{I,p}$, $1\le p\le\infty$, and $s\in \N$. Then
\begin{equation*}
%\label{NS}
   \sum_{[\beta]=s} \Vert D^\beta P\Vert_p\le c\,\omega_s(P,1)_p,
\end{equation*}
where $c$ does not depend on $P$.
\end{lem}

\begin{lem}\label{lem42}
  Let $P\in \mathcal{B}_{I,p}$, $1< p<\infty$, and $s>0$. Then
	$$
	c_1 \omega_s(P,1)_p\le\Vert (-\Delta)^{s/2} P\Vert_p\le c_2 \omega_s(P,1)_p,
	$$
where the constants $c_1$ and $c_2$ do not depend on $P$.
\end{lem}

\section{Main results}

\subsection{Main lemma}

Let $M\in \mathfrak{M}$, $\a\in \A_M$, and let  $\w\vp$ belong to $\mathcal{S}_{\a,p;M}'$. In what follows, we  understand
$\langle f,\w\vp_{jk}\rangle$ in the sense of Definition~\ref{def0}. 	
Thus, the quasi-projection  operators
$$
Q_j(f,\phi,\w\phi)=\sum_{k\in\zd}\langle f,\w\vp_{jk}\rangle\phi_{jk}
$$
are  defined for all $f\in \mathbb{B}_{p;M}^{\a(\cdot)}$.
By Lemmas~\ref{lem1} and~\ref{lem99}, we have that $\{\langle f,\w\vp_{jk}\rangle\}_k\in\ell_p$ and    $\{\langle f,\w\vp_{jk}\rangle\}_k\in {\rm c}_0$ if $p=\infty$. This, together with Lemmas~\ref{lemKK1+} and \ref{lemKK1}, implies that
the series $\sum_{k\in\zd}\langle f,\w\vp_{jk}\rangle\phi_{jk}$
converges unconditionally in $L_p$. Thus, the operator $Q_j(f,\phi,\w\phi)$ is well defined.

\smallskip

An analogue of the following lemma for the case $\phi\in {\cal L}_p$,
$\alpha(M)=|\det M|^N$, and $p<\infty$  can be found in~\cite{KS19}. In the general case, the proof is similar, but for completeness, we present it in detail.

\begin{lem}
\label{thKS}
Let $1\le p\le \infty$, $M\in \mathfrak{M}$,  $\delta\in (0,1]$, and $\nu\in\z_+$.  Suppose that $\phi\in {\cal L}_p$ or $\vp\in L_p$ is band-limited with $\h\vp \in \mathcal{M}_p$, and the functions
$P_\mu$, $\mu \in \Z_+$, are as in Definition~\ref{def0}.

\begin{itemize}
  \item[$(i)$] If $\a\in \A_M$, $\w\vp\in \mathcal{S}_{\a,p;M}'$, and $f\in \mathbb{B}_{p;M}^{\a(\cdot)}$,   then
\begin{equation}\label{KS000}
\begin{split}
\|f-Q_0(f,\vp,\w\vp)\|_p\le
\|P_\nu-Q_0(P_\nu,\vp,\w\vp)\|_p+
c\sum_{\mu=\nu}^\infty m^{\frac \mu p} \alpha (M^\mu) E_{\delta_p M^\mu}(f)_p.
\end{split}
\end{equation}

\item[$(ii)$] If  $\w\vp \in \mathcal{L}_{p'}$ and
$f\in L_p$, $p<\infty$, or
 $\w\vp\in \mathcal{S}_{\const,\infty;M}'$ and
$f\in {C}_0$, $p=\infty$, then
\begin{equation}\label{KS000_K}
\begin{split}
\|f-Q_0(f,\vp,\w\vp)\|_p\le
\|P_\nu-Q_0(P_\nu,\vp,\w\vp)\|_p+
c\,E_{\delta_p M^\nu}(f)_p.
\end{split}
\end{equation}
\end{itemize}
In the above two inequalities, the constant $c$ does not depend on $f$ and $\nu$.
\end{lem}
{\bf Proof.}
Obviously,
\begin{equation}\label{_1}
  \begin{split}
     \|f-Q_0(f,\vp,\w\vp)\|_p\le \|P_\nu-Q_0(P_\nu,\vp,\w\vp)\|_p+\|f-P_\nu\|_p+\|Q_0(f-P_\nu,\vp,\w\vp)\|_p.
  \end{split}
\end{equation}
%Thus, to prove the lemma, we need only to estimate $\|Q_0(f-P_\nu,\vp,\w\vp)\|_p$.
Suppose that conditions of item $(i)$ hold. Then, using
Lemmas~\ref{lemKK1+}, \ref{lemKK1}, and \ref{lem1},  we obtain
\begin{equation}\label{_2}
  \begin{split}
     \|Q_0(f-P_\nu,\vp,\w\vp)\|_p&\le c_1\|\{\langle f-P_\nu, \w\phi_{0k}\rangle\}_k\|_{\ell_p}\le c_2\sum_{\mu=\nu}^\infty\|\{\langle P_{\mu+1}-P_\mu, \w\phi_{0k}\rangle\}_k\|_{\ell_p}\\
&\le c_3\sum_{\mu=\nu}^\infty m^{\frac \mu p} \alpha (M^\mu) E_{\delta_p M^\mu}(f)_p.
  \end{split}
\end{equation}
Combining~\eqref{_1}, \eqref{_2}, and~\eqref{eT}, we get~\eqref{KS000}.

Similarly, under assumptions of item $(ii)$ in the case $p<\infty$, it follows from Lemmas~\ref{lemKK1+}, \ref{lemKK1}, and~\ref{lemNU} that
\begin{equation}\label{_3}
  \begin{split}
     \|Q_0(f-P_\nu,\vp,\w\vp)\|_p&\le c_4\|\{\langle f-P_\nu, \w\phi_{0k}\rangle\}_k\|_{\ell_p}\\
&\le c_4\Vert \w\vp \Vert_{\mathcal{L}_{p'}}\Vert f-P_\nu\Vert_p\le c_5E_{\delta_p M^\nu}(f)_p.
  \end{split}
\end{equation}
Thus, combining~\eqref{_1}, \eqref{_3} and~\eqref{eT}, we obtain~\eqref{KS000_K} for $p<\infty$.  In the case
$p=\infty$, using Lemma~\ref{lem99}
together with Lemma~\ref{lemKK1+},  we get
$$
 \|Q_0(f-P_\nu,\vp,\w\vp)\|_\infty\le c_6\|\{\langle f-P_\nu, \w\phi_{0k}\rangle\}_k\|_{\ell_\infty}\le
c_7E_{\delta_p M^\nu}(f)_\infty,
$$
which completes the proof of the lemma.~~$\Diamond$

\subsection{The case of weak compatibility of  $\phi$ and $\w\phi$}

In this subsection,
we give  error estimates for the quasi-projection operators associated with weakly compatible
 $\phi$ and  $\w\phi$.

\begin{theo}
\label{corMOD1'+}
Let $1\le p\le\infty$, $M\in \mathfrak{M}$, $\a\in\A_M$,  $s\in \N$, and $\delta\in (0,1]$.
Suppose that  $\w\phi\in \S'_{\a,p;M}$ and $\phi\in L_p$ satisfy the following conditions:
\begin{itemize}
  \item[$1)$] $\vp$ is band-limited with $\h\vp \in \mathcal{M}_p$ or $\phi\in{\cal L}_p$;
  \item[$2)$] the Strang-Fix condition of order $s$  holds for $\phi$;
  \item[$3)$] $\phi$ and ${\w\phi}$ are weakly  compatible of order~$s$;
  \item[$4)$] $\eta_\delta D^{\beta}\overline{\h\phi}{\h{\w\vp}} \in \mathcal{M}_{p}$ and $  \eta_\delta D^{\beta}\h{\phi}(\cdot+l)\in \mathcal{M}_p$ for all $\beta \in \Z_+^d$, $[\beta]=s$, and  $l\in\zd\setminus \{\nul\}$;
 \item[$5)$]  $ \sum_{l\ne\nul} \|\eta_\delta D^{\beta}\h{\phi}(\cdot+l)\|_{ \mathcal{M}_p}<\infty$ for all $\beta \in \Z_+^d$, $[\beta]=s$.
\end{itemize}
Then, for any $f\in \mathbb{B}_{p;M}^{\a(\cdot)}$, we have
\begin{equation}\label{110}
\begin{split}
\Vert f - Q_j(f,\phi,\w\phi) \Vert_p\le c\(\Omega_s(f, M^{-j})_p+m^{-\frac jp}\sum_{\nu=j}^\infty m^{\frac\nu p}\alpha(M^{\nu-j}) E_{M^\nu}(f)_p\).
\end{split}
\end{equation}
Moreover, if   $\w\vp \in \mathcal{L}_{p'}$ and
$f\in L_p$, $p<\infty$, or
 $\w\vp\in \mathcal{S}_{\const,\infty;M}'$ and
$f\in {C}_0$, $p=\infty$, then
\begin{equation}\label{110K}
  \Vert f - Q_j(f,\phi,\w\phi) \Vert_p\le c\,\Omega_s(f,M^{-j})_p.
\end{equation}
In the above inequalities, the constant $c$ does not depend on $f$ and $j$.
\end{theo}
{\bf Proof.} First we note that it suffices to prove~\eqref{110} and~\eqref{110K} for  $j=0$. Indeed,
	$$
\bigg\|f-\sum\limits_{k\in\zd}\langle f,\widetilde\phi_{jk}\rangle \phi_{jk}\bigg\|_p=
\bigg\|m^{-j/p}\big(f(M^{-j}\cdot)-\sum\limits_{k\in\zd}\langle f(M^{-j}\cdot),\widetilde\phi_{0k}\rangle \phi_{0k}\big)\bigg\|_p.
$$
Obviously, $m^{-j/p}f(M^{-j}\cdot)\in \mathbb{B}_{p;M}^{\a(\cdot)}$  whenever $f\in \mathbb{B}_{p;M}^{\a(\cdot)}$. We have also  that
$$
E_{ M^\nu}(m^{-j/p}f(M^{-j}\cdot))_p=E_{ M^{\nu+j}}(f)_p
$$
and
$$
\omega_s(f(M^{-j}\cdot), 1)_p= m^{j/p}\Omega_s(f, M^{-j})_p,
$$
which yields~\eqref{110} and~\eqref{110K} whenever these relations hold true for $j=0$.

Next, in view of Lemmas~\ref{thKS} and~\ref{lemJ}, to prove the theorem,  it suffices to show that
\be
\label{102}
\bigg\|P-\sum_{k\in\zd}\langle  P, \w\phi_{0k}\rangle\phi_{0k}\bigg\|_p \le c_1 \sum_{[\beta]=s}\|D^\beta P\|_p,
\ee
where the function $P$ is such that $P\in \mathcal{B}_{\d I,p}\cap L_2$  and
$
\|f-P\|_p\le c(d,p)E_{\d I}(f)_p
$
(remind that  $\delta_p=\delta$ for $p<\infty$ and $\delta_p=\delta/2$ for $p=\infty$, and
the function $P$ exists in view of Lemmas~\ref{lem2} and~\ref{lemJ2}).
Indeed, due to Lemmas~\ref{lemNS} and \ref{lemJ}, there holds
\begin{equation}\label{KKper1}
  \begin{split}
      \sum_{[\beta]=s} \Vert D^\b P \Vert_p\le &c_2\,\omega_s(P,1)_p\le c_2\Big(\omega_s(f,1)_p+E_{\delta_p I}(f)_p\Big)\le c_3\,\omega_s(f,1)_p.
   \end{split}
\end{equation}
Thus, combining~\eqref{KKper1} and~\eqref{102} with Lemma~\ref{thKS},  we obtain
$$
\Vert f - Q_0(f,\phi,\w\phi) \Vert_p\le c_4 \(\omega_s(f,1)_p+\sum_{\nu=0}^\infty m^{\frac\nu p}\alpha(M^{\nu}) E_{\delta_p M^\nu}(f)_p\).
$$
Since there exists $\nu_0=\nu(\d)\in \N$ such that $E_{\delta_p M^\nu}(f)_p\le E_{M^{\nu-\nu_0}}(f)_p$
and $\alpha(M^{\nu})\le c(\d,M)\alpha(M^{\nu-\nu_0})$
for all  $\nu>\nu_0$, applying Lemma~\ref{lemJ} and the inequality $\omega_s(f,\l)_p\le (1+\l)^s \omega_s(f,1)_p$ (see, e.g.,~\cite{KT20}) to the first~$\nu_0$ terms of the sum, we get~\eqref{110} for $j=0$.
Similarly, taking into account Lemmas~\ref{lemJ},  we derive~\eqref{110K}.
Thus, it remains to verifying inequality~\eqref{102}.

Set
$$
\Psi_0=1-\h{\vp}\overline{\h{\w\vp}}\quad\text{and}\quad \Psi_{l}=\h\vp(\cdot+l)\overline{\h{\w\vp}},\quad l\in \Z^d\setminus\{\nul\},
$$
and estimate $\|\Lambda_{\Psi_l}(P)\|_p$ for all $l\in\zd$.

Let $l\in \Z^d\setminus\{\nul\}$.  Using condition 2) and Taylor's formula, we have
\begin{equation*}
%\label{-11+}
  \h{\phi} (\xi+l) = \sum_{[\beta]=s} \frac{s}{\beta!}  \xi^{\beta} \int_0^1 (1-t)^{s-1} D^{\beta}\h{\phi}( t\xi+l) d t,
\end{equation*}
  which yields for $p<\infty$ that
  \begin{equation*}
    \begin{split}
              \|\Lambda_{\Psi_{l}}(P)\|_p^p&= \int_{\rd}\,dx\Big|\sum_{[\b]=s}\frac{s}{\b!} \int_0^1\,dt (1-t)^{s-1}
			\int_{\rd}\,d\xi\(D^\b\h\vp(t\xi+l)\xi^\b\h{P}(\xi)\overline{\h{\w\vp}(\xi)}\) e^{2\pi i(\xi,x)} \Big|^p
\\
&= \int_{\rd}\,dx\Big|\sum_{[\b]=s}\frac{s}{\b!(2\pi i)^{[\b]}}
			\int_0^1\,dt (1-t)^{s-1}
			\int_{\rd}\,d\xi \eta_\delta(t\xi) D^\b\h\vp(t\xi+l)\h\Theta_\beta({\xi})\, e^{2\pi i(\xi,x)} \Big|^p,
     \end{split}
  \end{equation*}	
where
$$
\Theta_\beta=\mathcal{F}^{-1}\(\h{D^\beta P }\,\overline{\h{\w\vp}}\).
$$		
Since $\eta_\delta(t\cdot) D^\b\h\vp(t\cdot+l)\in \mathcal{M}_p$ for every $t>0$ and
$\|\eta_\delta(t\cdot) D^\b\h\vp(t\cdot+l)\|_{\mathcal{M}_p}$ does not depend on $t$ (see property (iv) of Fourier multipliers), it follows from condition 4) and inequality~\eqref{DefS} that
\begin{equation}
\label{701}
\begin{split}
\|\Lambda_{\Psi_{l}}(P)\|_p&\le \sum_{[\b]=s}\sup_{t\in (0,1)}\Big\|
\mathcal{F}^{-1}\(\eta_\delta(t\cdot)D^\b\h\vp(t\cdot+l)\h{\Theta_\beta}\)\Big\|_p\\
&\le  \sum_{[\b]=s}\|\eta_\delta D^{\beta}\h{\phi}(\cdot+l)\|_{ \mathcal{M}_p}\|\Theta_\beta\|_p
\\
&=\sum_{[\b]=s}\|\eta_\delta D^{\beta}\h{\phi}(\cdot+l)\|_{ \mathcal{M}_p}\Vert \Lambda_{\F(\w\phi^-)}(D^\beta P) \Vert_p\\
&\le \alpha(\delta I) \sum_{[\b]=s} \|\eta_\delta D^{\beta}\h{\phi}(\cdot+l)\|_{ \mathcal{M}_p} \Vert D^\beta P \Vert_p.
 \end{split}
\end{equation}	
Similarly,
\begin{equation}
\label{702}
\begin{split}
\|\Lambda_{\Psi_{l}}(P)\|_\infty&\le
 \sum_{[\b]=s}\sup_{t\in (0,1)}\Big\|
\mathcal{F}^{-1}\(\eta_\delta(t\cdot)D^\b\h\vp(t\cdot+l)\h{\Theta_\beta}\)\Big\|_\infty
\\
&\le \sum_{[\b]=s}\|\eta_\delta D^{\beta}\h{\phi}(\cdot+l)\|_{ \mathcal{M}_\infty}\Vert \Lambda_{\F(\w\phi^-)}
			(D^\beta P) \Vert_\infty\\
&\le \alpha(\delta I) \sum_{[\b]=s} \|\eta_\delta D^{\beta}\h{\phi}(\cdot+l)\|_{ \mathcal{M}_\infty}\Vert D^\beta P \Vert_\infty.
 \end{split}
\end{equation}		
Combining relations~\eqref{701} and \eqref{702} with condition 5), we get	
\be
\label{711}
\sum_{l\ne\nul}\|\Lambda_{\Psi_{l}}(P)\|_p\le c_5 \sum_{[\b]=s} \Vert D^\beta P \Vert_p.
\ee

To estimate $\|\Lambda_{\Psi_0}( P)\|_p$, we note that by condition 4) and  Taylor's formula, there holds
\begin{equation*}
%\label{11+}
  \h{\phi} (\xi)\overline{\h{\w\phi}(\xi)}-1=\sum_{[\beta]=s}
	\frac{s}{\beta!}  \xi^{\beta}
	\int_0^1 (1-t)^{s-1} D^{\beta}\h{\phi}\overline{\h{\w\phi}}( t \xi) d t.
\end{equation*}
As above, we obtain for $p<\infty$ that
\begin{equation}
  \begin{split}
\label{714}
      & \|\Lambda_{\Psi_0}( P)\|_p=
			\\
			&\(\,\int_{\rd}\,dx\bigg|\sum_{[\b]=s}\frac{s}{\b!(2\pi i)^{[\beta]}}
 \int\limits_0^1\,dt (1-t)^{s-1}
			\int_{\rd}\,d\xi\,\eta_\delta(t\xi)D^\beta(\h{\phi}\overline{\h{\w\phi}})( t \xi)
			\h{D^\beta P}(\xi) e^{2\pi i(\xi,x)} \bigg|^p\)^{1/p}
			\\
&\le \sum_{[\b]=s}\sup_{t\in (0,1)}\left\|
\mathcal{F}^{-1}\(\eta_\delta(t\cdot)D^\beta(\h{\phi}\overline{\h{\w\phi}})( t\cdot)\h P\)\right\|_p
\le c_6 \sum_{[\beta]=s}\|D^\beta P \|_p,
   \end{split}
\end{equation}
and, similarly,
%\be
%\label{712}
$$
\|\Lambda_{\Psi_0}( P)\|_\infty\le \sum_{[\b]=s}\sup_{t\in (0,1)}\left\|
\mathcal{F}^{-1}\(\eta_\delta(t\cdot)D^\beta(\h{\phi}\overline{\h{\w\phi}})( t\cdot)\h P\)\right\|_\infty
\le c_7 \sum_{[\beta]=s}\|D^\beta P \|_\infty.
$$
%\ee

Next, we set $G(\xi):=\sum_{k\in\zd}\h P(\xi+k)\overline{\h{\w\phi}(\xi+k)}$ and prove that
\be
\label{710}
\sum_{k\in\zd}\langle  P, \w\phi_{0k}\rangle\phi_{0k}=\mathcal{F}^{-1}(G\h\phi).
\ee
First we consider the case $\phi\in {\cal L}_p$. Let $l\in \zd\setminus\{\nul\}$ and let $h_l$ denote the restriction of $\h\phi$ onto the set $\td-l$. Then
$$
\mathcal{F}^{-1}\(G h_l\)(x)=\int_{\td-l}
G(\xi) h_l(\xi)e^{2\pi i(x,\xi)}\, d\xi= \int_{\td}
\h P(\xi)\overline{\h{\w\phi}(\xi)} \h\phi(\xi+l) e^{2\pi i(x,\xi+l)}\, d\xi=\Lambda_{\Psi_l}(P)(x).
$$

Denote by $\Omega$ and $\Omega_N$  respectively the sum and the $N$-th cubic partial sum of $\sum_{l\in\zd}\mathcal{F}^{-1}\(G h_l\)$, which  converges in $L_p$ because of~\eqref{711}. Let us check that $\Omega=\mathcal{F}^{-1}(G\h\phi)$ in the distribution sense. Since $\h\phi$ is bounded,  the function $G h_l$ is in $L_2$, which yields that $\sum_{ \|l\|_\infty\le N}G h_l=\mathcal{F}\Omega_N$ almost everywhere. Hence, for every $g\in \mathcal{S}$, we have
$$
\langle \mathcal{F}^{-1}\(G\h\phi\)-\Omega_N, g \rangle
=\Big\langle G\h\phi-\sum_{ \|l\|_\infty\le N}G h_l, \h g \Big\rangle\too_{N\to\infty}0
$$
and, obviously,
$$
\langle \Omega-\Omega_N, g \rangle\too_{N\to\infty}0.
$$
Thus, the tempered distribution $\Omega$ coincides with $\mathcal{F}^{-1}(G\phi)$.  On the other hand,
using Lemma~1 from~\cite{KS} and cubic convergence of the Fourier series of $G$ in $L_2$-norm,
we have the equality
\be
\label{713}
\mathcal{F}\( \sum_{k\in\zd}\langle
P, \w\phi_{0k}\rangle\phi_{0k}\)=G\h\phi
\ee
in the distribution sense. Thus, the functions $\Omega$ and
$Q_0(f,\vp,\w\vp)$ coincide as distributions. But both $\Omega$ and $Q_0(f,\vp,\w\vp)$
are locally summable, hence, due to the du Bois-Reymond lemma,
these functions coincide almost everywhere, and so~\eqref{710} is proved.

Now consider the case of bandlimited $\phi$. Again~\eqref{713} holds true
in the distribution sense. Since $G$ is locally in $L_2$ and $\h\phi$ is  bounded and compactly supported, we have
$G\h\phi\in L$. Thus, again both the functions  $\mathcal{F}^{-1}(G\h\phi)$ and $\sum_{k\in\zd}\langle  P, \w\phi_{0k}\rangle\phi_{0k}$ are locally summable, which yields~\eqref{710}.

It follows from~\eqref{710} that
\begin{equation}\label{gl}
  \begin{split}
     &\sum_{k\in\zd}\langle  P, \w\phi_{0k}\rangle\phi_{0k}(x)
     =\int_{\rd }G(\xi)\h\phi(\xi) e^{2\pi i(\xi,x)}\,d\xi\\
     &=\sum_{l\in\zd}\int_{\td }G(\xi)\h\phi(\xi+l) e^{2\pi i(\xi+l,x)}\,d\xi
     =\sum_{l\in\zd}e^{2\pi i(l,x)}\int_{\td }\h P(\xi)\overline{\h{\w\phi}(\xi)}\h\phi(\xi+l) e^{2\pi i(\xi,x)}\,d\xi.
  \end{split}
\end{equation}
From this, taking into account that $P=\mathcal{F}^{-1}(\h P)$, we obtain
%\be
%\label{107}
$$
\bigg\|P-\sum_{k\in\zd}\langle  P, \w\phi_{0k}\rangle\phi_{0k}\bigg\|_p \le
\|\Lambda_{\Psi_0}( P)\|_p+\sum_{l\ne\nul} \|\Lambda_{\Psi_{l}}(P)\|_p,
$$
%\ee
which together with~\eqref{711} and \eqref{714}  yields~\eqref{102}. This completes the proof of the theorem.~~$\Diamond$

\begin{coro}
\label{coro1NN}
Let $1\le p\le \infty$,  $s\in \N$, $\delta\in (0,1]$, $M\in \mathfrak{M}$, $\a\in \A_M$, and let
$\phi\in {\cal L}_p$ and $\w\phi\in \S'_{\a,p;M}$.
Suppose that conditions 2) and 3) of Theorem~\ref{corMOD1'+} are satisfied and, additionally,

\begin{itemize}
  \item[$a)$]
 if $1<p<\infty$, we  suppose that  for some $k\in \n$, $k>  d|\frac1p-\frac12|$,
%\be
%\label{com2+NN}
$$
\overline{\h\phi}{\h{\w\vp}}\in
C^{s+k}(2\d \T^d),\quad
 \h\phi(\cdot+l)\in C^{s+k}(2\d \T^d)\quad\text{for all}\quad l\in\Z^d\setminus \{{\bf 0}\},
%\ee
$$
and
$$
\sum_{l\neq \bf{0}}\sup_{\xi\in 2\d \T^d}|D^\beta \h \vp(\xi+l)|^{1-\frac{d}{k}|\frac1p-\frac12|}<\infty\quad \text{for all}\quad \beta \in \Z_+^d,\quad [\beta]=s;
$$

\item[$b)$] if $p=1$ or $p=\infty$, we suppose that for some $k\in \N$, $k>\frac d2$,
%\be
%\label{com2+PrNN}
$$
\overline{\h\phi}{\h{\w\vp}}\in W_2^{s+k}(2\delta \T^d),
\quad \h\phi(\cdot+l) \in W_2^{s+k}(2\delta \T^d)\quad\text{for all}\quad l\in\Z^d\setminus \{{\bf 0}\},
$$
%\ee
and
$$
\sum_{l\neq \bf{0}}\Vert D^\beta \h \vp(\cdot+l)\Vert_{L_2(2\d \T^d)}^{1-\frac{d}{2k}}<\infty\quad \text{for all}\quad \beta \in \Z_+^d,\quad [\beta]=s.
$$
\end{itemize}
Then inequalities~\eqref{110} and~\eqref{110K} hold true.

For a band-limited function $\vp$, the above statement remains valid if the condition $\vp\in {\cal L}_p$ is replaced by the assumption that $\vp\in L_1$ in the case $p=1,\infty$ and $\vp=\F^{-1}(\chi_U\psi)$, where $U$ is compact and $\psi\in C^k(\R^d)$ if $1<p<\infty$.
\end{coro}

The proof of Corollary~\ref{coro1NN} easily follows from sufficient conditions for Fourier multipliers given
in~\cite[Corollaries~2 and 3]{Kol}. Let us compare this result with Theorem~10 in~\cite{KS19}, where the same estimates are obtained for the case $\alpha(M)=|\det M|^N$, but the proof is given without using Fourier multipliers.
A higher order of smoothness near the integer points is required for functions $\h\phi$ and $\h{\w\phi}$ in that theorem,
namely, differentiability  of order $s+d+1$ is assumed. However, the requirement on the decay for the functions $D^\beta\h\phi$
is less restrictive there. To provide the same decay for $D^\beta\h\phi$, we give another corollary based on Mikhlin's condition (see Section~3) for Fourier multipliers.

\begin{coro}
\label{coro1NN+}
Let $1< p< \infty$,  $s\in \N$, $M\in \mathfrak{M}$, $\a\in \mathcal{A}_M$, $\d\in (0,1]$, $k\in \n$, $k>\frac d2$, and let  $\w\phi\in \S'_{\a,p;M}$ and $\phi\in{\cal L}_p$.
Suppose that conditions 2) and 3) of Theorem~\ref{corMOD1'+} are satisfied. Additionally,
$$
\overline{\h\phi}{\h{\w\vp}}\in
C^{s+k}(2\d \T^d),\quad
 \h\phi(\cdot+l)\in C^{s+k}(2\d \T^d)\quad\text{for all}\quad l\in\Z^d\setminus \{{\bf 0}\},
$$
and
$$
\sum_{l\neq \bf{0}}\sup_{\xi\in 2\d \T^d}|D^\beta \h \vp(\xi+l)|<\infty\quad \text{for all}\quad \beta \in \Z_+^d,\quad [\beta]=s.
$$
Then inequalities~\eqref{110} and~\eqref{110K} hold true.

For a band-limited function $\vp$, the above statement remains valid if the condition $\vp\in {\cal L}_p$ is replaced by the assumption that $\vp=\F^{-1}(\chi_U\psi)$, where $U$ is compact and $\psi\in C^k(\R^d)$.
\end{coro}

\noindent\textbf{Example 1.} Let $d=2$, $p=\infty$,  $s=2$, $f\in {C}_0$, and let
$Q_j(f,\phi,\w\phi)$   be a  mixed sampling-Kantorovich  quasi-projection operator associated with
$\phi(x)=\frac1{16}\sinc^3(\frac{x_1}4)\sinc^3(\frac{x_2}4)$ and
$\w\phi(x)=\delta(x_1)\chi_{\T^1}(x_2)$, i.e.,
$$
Q_0(f,\vp,\w\vp)(x)=\sum_{k\in\z^2}\int_{k_2-1/2}^{k_2+1/2} f(k_1, t)\,dt\, \phi (x-k).
$$
It is easy to see that all assumptions of Theorem~\ref{corMOD1'+} for the case $p=\infty$ and
$\w\vp\in \S_{{\rm const},\infty;M}'$ are satisfied, which implies
$$
\|f-Q_j(f,\vp,\w\vp)\|_\infty\le c\,\Omega_2(f,M^{-j})_\infty.
$$

\medskip

%In the next example, we consider a reverse situation in some sense.
%
%\medskip

\noindent\textbf{Example 2.} Let $d=2$, $1<p<\infty$,   $s\in \N$, $f\in L_p$, and let $Q_j(f,\phi,\w\phi)$ be a  differential sampling expansion of the form
$$
Q_j(f,\vp,\w\vp)(x)=\sum_{k\in\z^2}\(f(M^{-j}k)+\frac{\partial^s}{\partial x_1^s}f(M^{-j}k)\)\sinc(M^jx-k),
$$
i.e., $\vp(x)=\sinc x$ and $\w\vp(x)=(I+\frac{\partial^s}{\partial x_1^s}) \d(x)$. We have $\w\vp\in \S_{\a,p;M}'$, where $\a(M)=m_1^s$ in the case  of the diagonal dilation matrix $M={\rm diag}\{m_1,m_2\}$. Thus, using Theorem~\ref{corMOD1'+}, it is not difficult to see that
$$
\Vert f-Q_j(f,\vp,\w\vp)\Vert_p\le c\(\Omega_s(f,M^{-j})_p+m_1^{-(s+\frac1p)j}m_2^{-\frac jp}\sum_{\nu=j}^\infty m_1^{(s+\frac1p)\nu}m_2^{\frac{\nu}p}E_{M^\nu}(f)_p\),
$$
where $f\in \mathbb{B}_{p;M}^{\a(\cdot)}$ and $c$ does not depend on $f$ and $j$.

\subsubsection*{Fractional smoothness and lower estimates}

We have the following generalization of Theorem~\ref{corMOD1'+} in terms of fractional moduli of smoothness.

\begin{theo}
\label{corMOD1+f++}
Let $1< p<\infty$, $s>0$, $\delta\in (0,1/2)$, $M\in \mathfrak{M}$, and  $\a\in \A_M$.
Suppose that $\w\phi\in \S'_{\a,p;M}$, $\phi\in L_p$ and  the following is satisfied:
\begin{itemize}
  \item[$1)$] $\supp\h\vp \subset\T^d$ and $\h\vp \in \mathcal{M}_p$;
  \item[$2)$] $\eta_\delta\frac{1-\overline{\h\phi}{\h{\w\vp}}}{|\cdot|^s} \in \mathcal{M}_p$. %and $\eta_\delta \frac{\h{\phi}(\cdot+l)}{|\cdot|^s} \in \mathcal{M}_p$ for any $l\in \Z^d\setminus \{\bf{0}\}$.
\end{itemize}
Then, for any $f\in \mathbb{B}_{p;M}^{\a(\cdot)}$, we have
\begin{equation*}
%\label{110++S}
\begin{split}
\Vert f - Q_j(f,\phi,\w\phi) \Vert_p\le c \(\Omega_s(f, M^{-j})_p+m^{-\frac jp}\sum_{\nu=j}^\infty m^{\frac\nu p}\alpha(M^{\nu-j}) E_{M^\nu}(f)_p\).
\end{split}
\end{equation*}
Moreover, if $\w\vp \in \mathcal{L}_{p'}$, then for any $f\in L_p$, we have
\begin{equation*}
%\label{110K++S}
  \Vert f - Q_j(f,\phi,\w\phi) \Vert_p\le c\, \Omega_s(f,M^{-j})_p.
\end{equation*}
In the above inequalities, the constant $c$ does not depend on $f$ and $j$.
\end{theo}
\textbf{Proof.}
%The proof is very similar to the proof of Theorem~\ref{corMOD1'+}.
Repeating the arguments of the proof of Theorem~\ref{corMOD1'+} for the case of bandlimited $\phi$, we see that it suffices to verify that
for any $P \in \mathcal{B}_{\d I,p}\cap L_2$ such that $\Vert f-P\Vert_p\le c(d,p)E_{\d I}(f)_p$, one has
\begin{equation}\label{eqNN1}
  \|\Lambda_{\Psi_0}(P)\|_p\le c_1 \omega_s(f,1)_p,
\end{equation}
where
$
\Psi_0=1-\h{\vp}\overline{\h{\w\vp}}.%,\quad  \Psi_{l}=\h\vp(\cdot+l)\overline{\h{\w\vp}}.
$

Using condition 2),  we derive
\begin{equation*}
  \begin{split}
%       \sum_{l\ne\nul}\|\Lambda_{\Psi_{l}}(P)\|_p&\le \a(\d I) \sum_{l\ne\nul}\|\mathcal{F}^{-1} (\h\vp(\cdot+l) \h P)\|_p\\
%       &\le c_3 \max_{l\ne\nul}\Big\|\mathcal{F}^{-1} \(\h\vp(\cdot+l)|\cdot|^{-s}\eta_\delta |\cdot|^s\h P\)\Big\|_p\le c_4 \Vert (-\Delta)^{s/2} P\Vert_p.
       \|\Lambda_{\Psi_{0}}(P)\|_p&=\Big\|\mathcal{F}^{-1} \((1-\h{\vp}\overline{\h{\w\vp}})|\cdot|^{-s}\eta_\delta |\cdot|^s\h P\)\Big\|_p
       %&\le c_3 \max_{l\ne\nul}\Big\|\mathcal{F}^{-1} \(\h\vp(\cdot+l)|\cdot|^{-s}\eta_\delta |\cdot|^s\h P\)\Big\|_p
       \le c_2 \Vert (-\Delta)^{s/2} P\Vert_p.
   \end{split}
\end{equation*}
Thus, to get~\eqref{eqNN1}, it remains to note that, due to Lemmas~\ref{lem42} and~\ref{lemJ}, we have
$$
\Vert (-\Delta)^{s/2} P\Vert_p\le c_3\, \omega_s(P,1)_p\le c_4(\|f-P\|_p+\omega_s(f,1)_p)\le c_5\,\omega_s(f,1)_p,
$$
which proves the theorem.~~$\Diamond$

\bigskip

In the next theorem, we obtain lower estimates for the $L_p$-error of approximation by the quasi-projection operators $Q_j(f,\phi,\w\phi)$. Note that such type of estimates are also called strong converse inequalities, see, e.g.,~\cite{DI}.

\begin{theo}\label{corMOD2conv}
Let $1< p<\infty$, $s>0$, $M\in \mathfrak{M}$,  and $\a\in \A_M$.
Suppose that $\w\phi\in \S'_{\a,p;M}$ and $\phi\in L_p$ satisfy the following conditions:
\begin{itemize}
  \item[$1)$] $\supp \h\vp\subset \T^d$  and $\h\vp \in \mathcal{M}_p$;
  \item[$2)$]  $\eta\frac{|\cdot|^s}{1-\overline{\h\phi}{\h{\w\vp}}} \in \mathcal{M}_p$.
\end{itemize}
Then, for any $f\in \mathbb{B}_{p;M}^{\a(\cdot)}$, we have
\begin{equation}\label{110inv}
\begin{split}
\Omega_s(f, M^{-j})_p\le c\,\Vert f - Q_j(f,\phi,\w\phi) \Vert_p+c\, m^{-\frac jp}\sum_{\nu=j}^\infty m^{\frac\nu p}\alpha(M^{\nu-j}) E_{M^\nu}(f)_p.
\end{split}
\end{equation}
Moreover, if $\w\vp \in \mathcal{L}_{p'}$, then for any $f\in L_p$, we have
\begin{equation}\label{110Kinv}
  \Omega_s(f,M^{-j})_p\le c\,\Vert f - Q_j(f,\phi,\w\phi) \Vert_p.
\end{equation}
In the above inequalities, the constant $c$ does not depend on $f$ and $j$.
\end{theo}
\textbf{Proof.} As in the proof of the previous theorems, it suffices to consider only the case $j=0$.
Let $P \in \mathcal{B}_{I,p}\cap L_2$ be such that $\Vert f-P\Vert_p\le c(d,p)E_{I}(f)_p$.
Due to the same arguments as in the proof of Theorem~\ref{corMOD1'+}, we have~\eqref{gl}, which takes now the following form
$$
P-\sum\limits_{k\in \Z^d}
\langle P, \w\phi_{0k} \rangle \phi_{0k}=
\mathcal{F}^{-1}\(\h P\(1-\h\phi\overline{\h{\w\phi}}\)\).
$$
Using this equality, Lemma~\ref{lem42}, and condition~2), we derive
\begin{equation*}
  \begin{split}
      \Omega_s(f,I)_p&=\omega_s(f,1)_p\le \omega_s(P,1)_p+c_1\|f-P\|_p
			\le c_2\(\Vert (-\Delta)^{s/2} P\Vert_p+E_I(f)_p\)
\\
&=c_2\(\left\| \mathcal{F}^{-1}\(\frac{\eta|\cdot|^s}{1-\h\phi\overline{\h{\w\phi}}}
\h P(1-\h\phi\overline{\h{\w\phi}})\) \right\|_p+E_I(f)_p\)
\\
&\le c_3\(\bigg\Vert P-\sum\limits_{k\in \Z^d}
\langle P, \w\phi_{0k} \rangle \phi_{0k}\bigg\Vert_p+E_{I}(f)_p\)
\\
&\le c_3\(\Vert f-Q_0(f,\vp,\w\vp)\Vert_p + \|P-f\|_p+\bigg\Vert \sum\limits_{k\in \Z^d}
 \langle P-f, \w\phi_{0k} \rangle \phi_{0k}\bigg\Vert_p+E_{I}(f)_p\)
\\
&\le c_4\(\Vert f-Q_0(f,\vp,\w\vp)\Vert_p+E_{I}(f)_p+\Vert Q_0(f-P,\vp,\w\vp)\Vert_p\)\\
&\le c_5\(\Vert f-Q_0(f,\vp,\w\vp)\Vert_p+E_{I}(f)_p\),
   \end{split}
\end{equation*}
where the last inequality follows from Lemmas~\ref{lemKK1+} and~\ref{lemNU}.
Thus, to prove \eqref{110Kinv}, it remains to note that in view of the inclusion  $\supp\mathcal{F}\(Q_0(f,\vp,\w\vp)\)\subset \supp\h\phi\subset \td$, we have
$
E_{I}(f)_p\le \|f-Q_0(f,\vp,\w\vp)\|_p.
$

Similarly, using~\eqref{_2} instead of Lemmas~\ref{lemKK1+} and~\ref{lemNU},
one can prove~\eqref{110inv}.~~$\Diamond$

\medskip

\begin{rem}
Note that the conditions on functions/distributions $\vp$ and $\w\vp$ in Theorems~\ref{corMOD1+f++} and~\ref{corMOD2conv} can be also given in terms of smoothness of $\h\vp$ and $\h{\w\vp}$, similarly to those given in Corollaries~\ref{coro1NN} and~\ref{coro1NN+}. For this, one can use the sufficient conditions for Fourier multipliers mentioned in Section~3 as well as some results of papers~\cite{LST} and~\cite{Kol}.
\end{rem}

\noindent\textbf{Example 3.}
Let $1<p<\infty$, $\vp(x)=\sinc(x):=\prod_{\nu=1}^d \frac{\sin(\pi x_\nu)}{\pi x_\nu}$, and $\w\vp(x)=\chi_{\T^d}(x)$ (the characteristic function of $\T^d$). Then all conditions of Theorems~\ref{corMOD1+f++}  and~\ref{corMOD2conv}  are satisfied and, therefore, for any  $f\in L_p$, we have
\begin{equation}\label{KKM}
 \bigg\Vert f- \sum_{k\in \Z^d} {m^j}\bigg(\int_{M^{-j}\T^d} f(M^{-j}k-t) dt\bigg) \sinc(M^j \cdot - k)\bigg\Vert_p\asymp\Omega_2(f,M^{-j})_p,
\end{equation}
where $\asymp$ is a two-sided inequality with constants independent of $f$ and $j$. Note that if we replace $\sinc x$ by $\sinc^2 x$, then the upper estimate in~\eqref{KKM} via the modulus $\Omega_2(f,M^{-j})_p$ holds for all $f\in L_p$, $1\le p\le \infty$. This can be easily verified using Theorem~\ref{corMOD1'+}, see also~\cite{KS3}.

\medskip

Similarly, using Theorems~\ref{corMOD1+f++}  and~\ref{corMOD2conv} and some basic properties of Fourier multipliers (see Section~3 and also~\cite{RS10} for some special multipliers), one can prove the following $L_p$ error estimates for approximation by quasi-projection operators generated by the Bochner-Riesz kernel of fractional order.

\medskip

\noindent\textbf{Example 4.}
Let $1< p<\infty$ and $\vp(x)=R_{s}^\gamma(x):=\mathcal{F}^{-1}\((1-|3\xi|^s)_+^\g\)(x)$, $s>0$, $\g>\frac{d-1}{2}$.
  \begin{itemize}
    \item[1)] If $\w\vp(x)=\d(x)$, then for any $f\in \mathbb{B}_{p;M}^1$, we have
\begin{equation*}\label{ee61}
\begin{split}
c_1\Omega_s(f,M^{-j})_p\le \bigg\Vert f- {m^j}\sum_{k\in \Z^d} &f(M^{-j}k){R}_{s}^\gamma (M^j \cdot - k)\bigg\Vert_p\le c_2 m^{-\frac jp}\sum_{\nu=j}^\infty m^{\frac\nu p}\Omega_s(f,M^{-\nu})_p,
\end{split}
\end{equation*}
where $c_1$ and $c_2$ are some positive constants independent of $f$ and $j$
    \item[2)] If $\w\vp(x)=\chi_{\T^d}(x)$, then  for any $f\in  L_p$ and $s\in (0,2]$, we have
\begin{equation*}\label{ee62}
\begin{split}
\bigg\Vert f- {m^j}\sum_{k\in \Z^d} \bigg(\int_{M^{-j}\T^d} f(M^{-j}k-t) dt\bigg) &{R}_{s}^\gamma (M^j \cdot - k)\bigg\Vert_p\asymp\Omega_s(f,M^{-j})_p,
\end{split}
\end{equation*}
where $\asymp$ is a two-sided inequality with positive constants independent of $f$ and $j$.
  \end{itemize}

\subsection{The case of strict compatibility for  $\phi$ and $\w\phi$}

\begin{theo}\label{cor1}
Let $1\le p\le\infty$, $\d\in (0,1]$, $M\in \mathfrak{M}$, and $\a\in \A_M$.
Suppose that $\w\phi\in \mathcal{S}'_{\a,p;M}$ and $\phi\in L_p$ satisfy the following conditions:
\begin{itemize}
  \item[$1)$] $\supp \h\vp\subset \T^d$ and $\h\vp \in \mathcal{M}_p$;
  \item[$2)$] $\phi$ and ${\w\phi}$ are strictly  compatible with respect to $\d$.
\end{itemize}
Then, for any $f\in \mathbb{B}_{p;M}^{\a(\cdot)}$, we have
\begin{equation}%\label{KS000+}
\begin{split}
\label{703}
     \Vert f - Q_j(f,\phi,\w\phi) \Vert_p\le c\,m^{-\frac jp}\sum_{\nu=j}^\infty m^{\frac\nu p}\alpha(M^{\nu-j})
E_{\d_p M^\nu}(f)_p.
\end{split}
\end{equation}
Moreover, if $\w\vp \in \mathcal{L}_{p'}$ and
$f\in L_p$, $p<\infty$, or $\w\vp\in \S_{{\rm const},\infty;M}'$ and $f\in {C}_0$, $p=\infty$,  then

\begin{equation}\label{KS000KLLL}
  \Vert f - Q_j(f,\phi,\w\phi) \Vert_p\le c\, E_{\d_p M^j}(f)_p.
\end{equation}
In the above inequalities, the constant $c$ does not depend on $f$ and $j$.
\end{theo}
{\bf Proof.}
As above, it suffices to consider only the case $j=0$. Repeating the arguments of the proof of Theorem~\ref{corMOD1'+}, we obtain
from~\eqref{gl} that
$$
P_0-\sum\limits_{k\in\zd}\langle P_0,\widetilde\phi_{0k}\rangle \phi_{0k}=0,
$$
Thus, applying Lemma~\ref{thKS}, we prove both the statements  of the theorem.~~$\Diamond$

\begin{rem}
In the case $p=\infty$, $\d_p$ in estimates~\eqref{703} and~\eqref{KS000KLLL} can be replaced by $\rho \delta$, where $\rho \in (0,1)$.
\end{rem}

\noindent\textbf{Example 5.} If $\w\vp(x)=\chi_{\T^d}(x)$ and
$\vp(x)=\mathcal{F}^{-1}\(\frac{\chi_{\T^d}(\xi)}{\sinc(\xi)}\)(x)$,
then Theorem~\ref{cor1} provides the following estimate for the corresponding Kantorovich-type operator
\begin{equation*}
 \bigg\Vert f- \sum_{k\in \Z^d} {m^j}\bigg(\int_{M^{-j} \T^d} f(M^{-j}k-t) dt\bigg) \vp(M^j\cdot-k)\bigg\Vert_p\le c\,E_{M^j}(f)_p,
\end{equation*}
where $f\in L_p$, $1<p<\infty$, and $c$ does not depend on $f$ and $j$. Note that the corresponding estimate in the case of Kotelnikov operators (the case $\w\vp(x)=\d(x)$ and $\vp(x)=\sinc(x)$) has the following form:
$$
 \bigg\Vert f- m^j\sum_{k\in \Z^d} f(M^{-j}k) {\rm \sinc}(M^j\cdot-k)\bigg\Vert_p\le c\, m^{-\frac jp}\sum_{\nu=j}^\infty m^{\frac{\nu}p} E_{M^\nu}(f)_p,\quad f\in \mathbb{B}_{p;M}^{1}.
$$

\subsubsection*{Whittaker--Nyquist--Kotelnikov--Shannon-type theorems}

One can easily see that the right hand side of~\eqref{703} is identically zero if $\supp \h f\subset \delta_p M^*\td$ and the matrix $M$ is such that
$M\td\subset \td$. This leads to the following counterpart of the classical Kotelnikov
formula
$$
f=Q_1(f,\phi,\w\phi)\quad a.e.
$$

The next two theorems provide results of this type under significantly milder conditions.

\begin{theo}\label{koteln1}
Let $M$ be a non-degenerate matrix and $\delta\in (0,1]$. Suppose that
\begin{itemize}
  \item[$1)$] $\supp\h\vp\subset\T^d$ and $\h\phi\in L_\infty$;
	\item[$2)$]  $\w\phi\in \mathcal{S}'$ and  $\h{\w\phi}$ is bounded on $\delta \td$;
  \item[$3)$] $\phi$ and ${\w\phi}$ are strictly  compatible with  respect to $\delta$.
\end{itemize}
If a function $f$  is such that  $\supp \h f\subset\delta M^*\td$ and $\h f\in L_q$, $q>1$,  then
\be
\label{Kot1}
f(x)=\lim_{n\to\infty}\sum_{\|k\|_\infty\le n} \langle \widehat{f}, \h{\w\phi_{1k}}\rangle \phi_{1k}(x)\quad \text{for almost all }\quad x\in\rd.
\ee
\end{theo}
{\bf Proof.} First let $M=I$.
Set $G(\xi):=\sum_{k\in\zd}\h f(\xi+k)\overline{\h{\w\phi}(\xi+k)}$.
Since $G\in L_q(\td)$, its
Fourier series   is cubic convergent to $G$  in $L_q(\td)$,
i.e.,  $\|G-G_n\|_{L_q(\td)} \to 0$, where $G_n$ is the $n$-th cubic partial Fourier sum and
$$
G_n(\xi)\h\phi(\xi)=\sum\limits_{\|k\|_{\infty}\le n}
\h G(k)e^{2\pi i(k,\xi)}
\h\phi(\xi)= \sum\limits_{\|k\|_{\infty}\le n}\langle \h f,\h{{\w\phi}_{0k}}
\rangle\widehat{\phi_{0k}}(\xi)=:H_n(\xi).
$$
Obviously, $H_n\in L_q$. This, together with condition 1), yields
\begin{equation*}
  \begin{split}
      \|H_m-H_n\|^q_q &= \int_{\rd}|(G_m(\xi) - G_n(\xi))\h\phi(\xi)|^q\,d\xi\\
                    &=\int_{\td}|(G_m(\xi) - G_n(\xi))\h\phi(\xi)|^q\,d\xi\le \|\h\phi\|^q_\infty\|G_m - G_n\|_{L_q(\td)}^q.%\quad n,m\to\infty
   \end{split}
\end{equation*}
Thus, the sequence $\{H_n\}_n$ converges in $L_q$. Without loss of generality, we can suppose that $q\le2$. By the Hausdorff-Young inequality, $\mathcal{F}^{-1}H_n$ and $\mathcal{F}^{-1}H_m$ are in $L_p$, where $p=\frac q{q-1}$, and
$$
\|\mathcal{F}^{-1} H_m - \mathcal{F}^{-1}H_n\|_p \le \|H_m-H_n\|_q \to 0\quad \text{as} \quad n, m \to \infty.
$$
It follows  that the series
$\sum_{k\in\zd} \langle \h f,\h{\w\phi_{0k}}\rangle \phi_{0k}$
is  cubic convergent   in $L_p$ and its sum  is in~$L_p$.
Again by the Hausdorff-Young inequality,
\begin{equation*}
  \begin{split}
     \Big\|f-\lim_{n\to\infty}&\sum_{\|k\|_\infty\le n}\langle \widehat{f}, \h{\w\phi_{0k}}\rangle\phi(\cdot+k)\Big\|_p\le\|\h f-G\h\phi\|_q\\
     &\le
\|\h f(1-\h\phi\overline{\h{\w\phi}})\|_q+\Big\|\h\phi\sum_{k\ne\nul}
\h f(\cdot+k)\overline{\h{\w\phi}(\cdot+k)}\Big\|_q:= I_1+I_2.
   \end{split}
\end{equation*}
Using condition 3)  and taking into account that
$ \supp \h f\subset \delta \td$, we obtain that $I_1=0$. At the same time, we have that
$I_2=0$ because both the functions $\phi$ and $f$ are band-limited to $\td$.
This yields~\eqref{Kot1} for the case $M=I_d$.

Now let $M$ be an arbitrary non-degenerate matrix. Setting $g:=f(M^{-1}\cdot)$,
we have $\supp \h g\subset \delta\td$. Hence, equality~\eqref{Kot1} holds for $g$ and
$$
f(x)=g(Mx)=\lim_{n\to\infty}\sum_{\|k\|_\infty\le n}
\langle \widehat{g}, \h{\w\phi_{0k}}\rangle \phi_{0k}(Mx) \quad \text{for almost all }\quad x\in\rd.
$$
Finally, after a suitable change of variable in the inner products,
we get~\eqref{Kot1}.~~$\Diamond$

\medskip

There are two drawbacks in the latter theorem (as well as in the Kotelnikov-type formula
extracted from Theorem~\ref{cor1}). First, the Fourier transform of $f$
is assumed to be in $L_q$, $q>1$, and second,  the Kotelnikov-type equality holds
only at almost all points. Under an additional restriction on $\phi$, these drawbacks
are avoided in the following statement.

\begin{theo}\label{koteln2}
Let $M$ be a non-degenerate matrix and $\delta\in (0,1]$.   Suppose that
\begin{itemize}
  \item[$1)$] $\supp\h\vp\subset\td$ and $\h\phi\in L_\infty$;
  \item[$2)$] $\Big|\frac{\partial\h\phi}{\partial x_l}(\xi)\Big|\le B$ for all $\xi\in\rd$ and $l=1,\dots,d$;
	\item[$3)$] $\w\phi\in \S'$ and   $\h{\w\phi}$ is bounded on $\delta\td$;
  \item[$4)$] $\phi$ and ${\w\phi}$ are strictly  compatible with parameter $\delta$ .
\end{itemize}
If a function $f$ is such that  $\supp \h f\subset\delta M^*\td$  and $\h f\in L_1$, then
$$
f(x)=\lim_{n\to\infty}\sum_{\|k\|_\infty\le n}
\langle \widehat{f}, \h{\w\phi_{1k}}\rangle
\phi_{1k}(x)\quad\text{for all}\quad x\in \rd.
$$
\end{theo}
{\bf Proof.} First let $M=I$. Set
$\Theta_x(\xi):=\sum_{s\in\,\zd}\h\phi(\xi+s) e^{2\pi i(x, \xi+s)}$.
By the Poisson summation formula, $\Theta_x$ is a summable $1$-periodic (with respect to
each variable) function and its $n$-th Fourier coefficient is
$$
\h{\Theta_x}(k)=
\int_{\rd}
\h\phi(\xi)e^{-2\pi i(k-x,\xi)}\,d\xi=
\h{\h\phi}(k-x)=\phi(x-k).
$$
Since  $\Theta_x$ is a bounded function,
its Fourier series cubic
converges  almost everywhere.
Let us check that the cubic partial
Fourier sums are uniformly bounded in $L_\infty$-norm. Set
$$
S_n(\Theta_x, \xi):=\sum_{\|k\|_\infty\le n}
\h{\Theta_x}(k)e^{2\pi i (k,\xi)}.
$$
Using the Lebesgue inequality and the Jackson type inequality  for the rectangular best approximations of periodic functions (see~\cite[Sec.~5.3.1]{Timan}), we have
%\be
\begin{equation*}
  \begin{split}\label{901}
      \|S_n(\Theta_x, \cdot)\|_\infty&\le \|\Theta_x\|_\infty+\|\Theta_x-S_n(\Theta_x, \cdot)\|_\infty\\
      &\le\|\Theta_x\|_\infty + c(d)\log^d (n+1)\,\sum_{\nu=1}^d\omega_1^{(\nu)}\Big(\Theta_x,  \frac 1n\Big)_\infty,
   \end{split}
\end{equation*}
%\ee
where  $\omega_1^{(\nu)}\big(g, h\big)_\infty=\sup_{|t|\le h}\Vert \D_{{\rm e}_\nu t}^1 g\Vert_\infty$ and $\{{\rm e}_\nu\}_{\nu=0}^d$ is the standard basis in $\R^d$.
Since the function $\Theta_x$ and all its partial derivatives are bounded (uniformly with respect to $x$), there exists a constant $c_1$ such that
$$
\sum_{\nu=1}^d\omega_1^{(\nu)}\Big(\Theta_x, \frac 1n\Big)_\infty
\le \frac {c_1}n \quad \text{for all} \quad x\in \rd,
$$
which together with~\eqref{901} implies the required boundedness.
Using this with Lebesgue's dominated convergence theorem and taking into account
that  both the functions $\phi$ and $f$ are band-limited to $\td$, we derive
\begin{equation*}
  \begin{split}
      \lim_{n\to\infty}&\sum_{\|k\|_\infty\le n}\langle \widehat{f}, \h{\w\phi_{0k}}\rangle\phi(x+k)=
\lim_{n\to\infty}
\int_{\rd}\sum_{\|k\|_\infty\le n}\phi(x+k)e^{-2\pi i (k,\xi)}
\overline{\h{\w\phi}}(\xi)\h f(\xi)\,d\xi
\\
&=\int_{\rd}   \lim_{n\to\infty}  \sum_{\|k\|_\infty\le n}   \phi(x-k)e^{2\pi i (k,\xi)}
\overline{\h{\w\phi}}(\xi)\h f(\xi)\,d\xi
\\
&=\int_{\rd}
e^{2\pi i(x, \xi)} \sum_{s\in\,\zd}
\h\phi(\xi+s)e^{2\pi i(x, s)}\overline{\h{\w\phi}}(\xi)\h f(\xi)\,d\xi=
\int_{\td}\h\phi(\xi)\overline{\h{\w\phi}}(\xi)\h f(\xi)e^{2\pi i(x, \xi)}\,d\xi.
   \end{split}
\end{equation*}
Since $\supp\h f\subset \delta\td$, it follows from condition 4)  that
$$
f(x)-\lim_{n\to\infty}\sum_{\|k\|_\infty\le n}\langle \widehat{f}, \h{\w\phi_{0k}}\rangle\phi(x+k)  =\int_{\delta\td}(1-\h\phi(\xi)\overline{\h\phi}(\xi))\h f(\xi)e^{2\pi i(x, \xi)}\,d\xi=0
$$
for every $x\in\rd$.

So, the theorem is proved for the case $M=I$. For the general case, it remains
to repeat the arguments at the end of the proof of Theorem~\ref{koteln1}.
~~$\Diamond$

\begin{rem}
Condition 2)  in Theorem~\ref{koteln2} can be replaced by the assumption that $\h\phi\in{\rm Lip}\,\alpha$, $\alpha>0$, with respect to each variable.
In the case $d=1$, condition~2)  can be also replaced by
the requirement of bounded variation for $\h\phi$.
\end{rem}

\noindent {\bf Acknowledgments} This research was supported by Volkswagen Foundation in framework
of the project "From Modeling and Analysis to Approximation".
The first author was partially supported by the DFG project KO 5804/1-1.

%\end{document}

\end{document}